\documentclass{amsart}
\usepackage{amsthm,amsfonts,amssymb,amsmath,amscd
 ,newclude
}
\usepackage{textcomp}
\usepackage{appendix}
\usepackage{xcolor}
\usepackage{longtable}
\usepackage{url}
\usepackage[centerlast]{subfigure}
\usepackage{esint}
\usepackage[all,2cell, cmtip]{xy} \UseAllTwocells

\begin{document}

\title[CoCartesian fibrations and homotopy colimits]{CoCartesian fibrations and homotopy colimits}
\author[A. Sharma]{Amit Sharma}
\email{asharm24@kent.edu}
\address {Department of mathematical sciences\\ Kent State University\\
  Kent,OH}

\date{\today}
%

\newcommand{\CONT}{\noindent}
\newcommand{\FIG}{Fig.\ }
\newcommand{\FIGS}{Figs.\ }
\newcommand{\SEC}{Sec.\ }
\newcommand{\SECS}{Secs.\ }
\newcommand{\TAB}{Table }
\newcommand{\TABS}{Tables }
\newcommand{\EQ}{Eq.\ }
\newcommand{\EQS}{Eqs.\ }
\newcommand{\APP}{Appendix }
\newcommand{\APPS}{Appendices }
\newcommand{\CHP}{Chapter }
\newcommand{\CHPS}{Chapters }

\newcommand{\OFF}{\emph{G2off}~}
\newcommand{\TOO}{\emph{G2Too}~}
\newcommand{\CatS}{Cat_{\bigS}}
\newcommand{\PicS}{{\underline{\pic}}^{\oplus}}
\newcommand{\HPicS}{{Hom^{\oplus}_{\pic}}}

\newtheorem{thm}{Theorem}[section]
\newtheorem{lem}[thm]{Lemma}
\newtheorem{conj}[thm]{Conjecture}
\newtheorem{coro}[thm]{Corollary}
\newtheorem{prop}[thm]{Proposition}

\theoremstyle{definition}
\newtheorem{df}[thm]{Definition}
\newtheorem{nota}[thm]{Notation}

\newtheorem{ex}[thm]{Example}
\newtheorem{exs}[thm]{Examples}

\theoremstyle{remark}
\newtheorem*{note}{Note}
\newtheorem{rem}{Remark}
\newtheorem{ack}{Acknowledgments}

\newcommand{\ChI}{{\textit{\v C}}\textit{ech}}
\newcommand{\Ch}{{\v C}ech}

\newcommand{\ChZG}{hermitian line $0$-gerbe}
\renewcommand{\theack}{$\! \! \!$}

\newcommand{\ChG}{flat hermitian line $1$-gerbe}
\newcommand{\ChC}{hermitian line $1$-cocycle}
\newcommand{\ChGG}{flat hermitian line $2$-gerbe}
\newcommand{\ChCC}{hermitian line $2$-cocycle}
\newcommand{\id}{id}
\newcommand{\LC}{\mathfrak{C}}
\newcommand{\Coker}{Coker}
\newcommand{\Com}{Com}
\newcommand{\Hom}{Hom}
\newcommand{\Mor}{Mor}
\newcommand{\Map}{Map}
\newcommand{\alg}{alg}
\newcommand{\an}{an}
\newcommand{\Ker}{Ker}
\newcommand{\Ob}{Ob}
\newcommand{\Proj}{\mathbf{Proj}}
\newcommand{\topo}{\mathbf{Top}}
\newcommand{\kan}{\mathcal{K}}
\newcommand{\pkan}{\mathcal{K}_\bullet}
\newcommand{\Kan}{\mathbf{Kan}}
\newcommand{\pKan}{\mathbf{Kan}_\bullet}
\newcommand{\QCat}{\mathbf{QCat}}
\newcommand{\gp}{\mathcal{A}_\infty}
\newcommand{\mdl}{\mathcal{M}\textit{odel}}
\newcommand{\sSets}{\mathbf{\S}}
\newcommand{\sSSets}{\sSets^{(2)}}
\newcommand{\Sets}{\mathbf{Sets}}
\newcommand{\sSetsM}{\mathbf{\S}^+}
\newcommand{\sSetsQ}{(\mathbf{\sSets, Q})}
\newcommand{\sSetsMQ}{(\mathbf{\sSetsM, Q})}
\newcommand{\sSetsK}{(\mathbf{\sSets, \Kan})}
\newcommand{\pSSets}{\mathbf{\sSets}_\bullet}
\newcommand{\pSSetsK}{(\mathbf{\sSets}_\bullet, \Kan)}
\newcommand{\pSSetsQ}{(\mathbf{\sSets_\bullet, Q})}
\newcommand{\cyl}{\mathbf{Cyl}}
\newcommand{\lin}{\mathcal{L}_\infty}
\newcommand{\Vect}{\mathbf{Vect}}
\newcommand{\Aut}{Aut}
\newcommand{\pic}{\mathcal{P}\textit{ic}}
\newcommand{\Dlin}{\pic}
\newcommand{\bigS}{\mathbf{S}}
\newcommand{\bigA}{\mathbf{A}}
\newcommand{\bhom}{\mathbf{hom}}
\newcommand{\bhomK}{\mathbf{hom}({\textit{K}}^+,\textit{-})}
\newcommand{\Bhom}{\mathbf{Hom}}
\newcommand{\bhomk}{\mathbf{hom}^{{\textit{k}}^+}}
\newcommand{\Dlino}{\pic^{\textit{op}}}
\newcommand{\lino}{\mathcal{L}^{\textit{op}}_\infty}
\newcommand{\lind}{\mathcal{L}^\delta_\infty}
\newcommand{\linK}{\mathcal{L}_\infty(\kan)}
\newcommand{\linC}{\mathcal{L}_\infty\text{-category}}
\newcommand{\linCs}{\mathcal{L}_\infty\text{-categories}}
\newcommand{\ainCs}{\text{additive} \ \infty-\text{categories}}
\newcommand{\ainC}{\text{additive} \ \infty-\text{category}}
\newcommand{\inC}{\infty\text{-category}}
\newcommand{\inCs}{\infty\text{-categories}}
\newcommand{\gS}{{\Gamma}\text{-space}}
\newcommand{\gSet}{{\Gamma}\text{-set}}
\newcommand{\ggS}{\Gamma \times \Gamma\text{-space}}
\newcommand{\gSs}{\Gamma\text{-spaces}}
\newcommand{\gSets}{\Gamma\text{-sets}}
\newcommand{\ggSs}{\Gamma \times \Gamma\text{-spaces}}
\newcommand{\gO}{\Gamma-\text{object}}
\newcommand{\gSCat}{{\Gamma}\text{-space category}}
\newcommand{\pss}{\mathbf{S}_\bullet}
\newcommand{\gSC}{{{{\Gamma}}\mathcal{S}}}
\newcommand{\gSCM}{{{{\Gamma}}\mathcal{S}^+}}
\newcommand{\pGSC}{{{{\Gamma}}\mathcal{S}}_\bullet}
\newcommand{\pGSCStr}{{{{\Gamma}}\mathcal{S}}_\bullet^{\textit{str}}}
\newcommand{\ggSC}{{\Gamma\Gamma\mathcal{S}}}
\newcommand{\gSD}{\mathbf{D}(\gSC^{\textit{f}})}
\newcommand{\sCat}{\mathbf{sCat}}
\newcommand{\pSCat}{\mathbf{sCat}_\bullet}
\newcommand{\gSetCat}{{{{\Gamma}}\mathcal{S}\textit{et}}}
\newcommand{\Dhom}{\mathbf{R}Hom_{\pic}}
\newcommand{\gop}{\Gamma^{\textit{op}}}
\newcommand{\fU}{\mathbf{U}}
\newcommand{\cDN}{\underset{\mathbf{D}[\textit{n}^+]}{\circ}}
\newcommand{\cDK}{\underset{\mathbf{D}[\textit{k}^+]}{\circ}}
\newcommand{\cDL}{\underset{\mathbf{D}[\textit{l}^+]}{\circ}}
\newcommand{\cD}{\underset{\gSD}{\circ}}
\newcommand{\cDT}{\underset{\gSD}{\widetilde{\circ}}}
\newcommand{\ppsSets}{\sSets_{\bullet, \bullet}}
\newcommand{\gdHom}{\underline{Hom}_{\gSD}}
\newcommand{\HomU}{\underline{Hom}}
\newcommand{\ominf}{\Omega_\infty}
\newcommand{\ev}{ev}
\newcommand{\cu}{C(X;\mathfrak{U}_I)}
\newcommand{\Sing}{Sing}
\newcommand{\AlgEin}{\A\textit{lg}_{\E_\infty}}
\newcommand{\SFunc}[2]{\mathbf{SFunc}({#1} ; {#2})}
\newcommand{\unit}[1]{\mathrm{1}_{#1}}
\newcommand{\liminj}{\varinjlim}
\newcommand{\limproj}{\varprojlim}
\newcommand{\HMapC}[3]{\mathcal{M}\textit{ap}^{\textit{h}}_{#3}(#1, #2)}
\newcommand{\tensPGSR}[2]{#1 \underset{\gSR}\wedge #2}
\newcommand{\pTensP}[3]{#1 \underset{#3}\wedge #2}
\newcommand{\MGCat}[2]{\underline{\map}_{\gSC}({#1},{ #2})}
\newcommand{\MGBoxCat}[2]{\underline{\map}_{\gSC}^{\Box}({#1},{ #2})}
\newcommand{\TensPFunc}[1]{- \underset{#1} \otimes -}
\newcommand{\TensP}[3]{#1 \underset{#3}\otimes #2}
\newcommand{\MapC}[3]{\mathcal{M}\textit{ap}_{#3}(#1, #2)}
\newcommand{\bHom}[3]{{#2}^{#1}}
\newcommand{\gn}[1]{\Gamma^{#1}}
\newcommand{\gnk}[2]{\Gamma^{#1}({#2}^+)}
\newcommand{\gnf}[2]{\Gamma^{#1}({#2})}
\newcommand{\ggn}[1]{\Gamma\Gamma^{#1}}
\newcommand{\Nat}{\mathbb{N}}
\newcommand{\partition}[2]{\delta^{#1}_{#2}}
\newcommand{\inclusion}[2]{\iota^{#1}_{#2}}
\newcommand{\EinQC}{\text{coherently commutative monoidal quasi-category}} 
\newcommand{\EinQCs}{\text{coherently commutative monoidal quasi-categories}}
\newcommand{\pHomCat}[2]{[#1,#2]_{\bullet}}
\newcommand{\CatHom}[3]{[#1,#2]^{#3}}
\newcommand{\pCatHom}[3]{[#1,#2]_\bullet^{#3}}
\newcommand{\EinC}{\text{coherently commutative monoidal category}}
\newcommand{\EinCs}{\text{coherently commutative monoidal categories}}
\newcommand{\EinLO}{E_\infty{\text{- local object}}}
\newcommand{\EinSLO}{\E_\infty\S{\text{- local object}}}
\newcommand{\Ein}{E_\infty}
\newcommand{\EinS}{E_\infty{\text{- space}}}
\newcommand{\EinSs}{E_\infty{\text{- spaces}}}
\newcommand{\PCat}{\mathbf{Perm}}
\newcommand{\nor}[1]{{#1}^\textit{nor}}
\newcommand{\pSSetsHom}[3]{[#1,#2]_\bullet^{#3}}
\newcommand{\PNat}{\overline{\L}}
\newcommand{\PStr}{\L}
\newcommand{\Gn}[1]{\Gamma[#1]}
\newcommand{\GIH}{\Gamma\textit{H}_{\textit{in}}}
\newcommand{\QStr}[1]{\L_\bullet(\ud{#1})}
\newcommand{\QStrF}{\L_\bullet}
\newcommand{\Kbar}{\overline{\K}}
\newcommand{\gPerm}{{\Gamma\PCat}}
\newcommand{\gCat}{{\Gamma\Cat}}
\newcommand{\MapS}[3]{\map_{#3}(#1, #2)}
\newcommand{\sSetsMG}{\sSetsM / N(\gop)}
\newcommand{\sSetsMGen}[1]{\sSetsM / N(#1)}
\newcommand{\pF}{\mathfrak{F}_\bullet^+(\gop)}
\newcommand{\pN}{{\textit{N}}_\bullet^+(\gop)}
\newcommand{\pFX}[1]{\mathfrak{F}_{#1}^+(\gop)}
\newcommand{\pNX}[1]{{\textit{N}}_{#1}^+(\gop)}
\newcommand{\nGop}{N(\gop)}
\newcommand{\sSetsMGSM}{(\sSetsM/ N(\gop), \otimes)}
\newcommand{\sSetsGen}[1]{\sSets/ #1}
\newcommand{\ovCatGen}[2]{#1/ #2}
\newcommand{\coMdl}[1]{(\sSets/ #1, \mathbf{L})}
\newcommand{\gCLM}[1]{\gCL{#1}^+}
\newcommand{\gCL}[1]{h_!}
\newcommand{\pFGen}[2]{\mathfrak{F}_{#1}^+(#2)}

\def\Pic{\mathbf{2}\mathcal P\textit{ic}}
\def\nc{\mathbb C}

\def\Z{\mathbb Z}
\def\P{\mathbb P}
\def\J{\mathcal J}
\def\I{\mathcal I}
\def\nC{\mathbb C}
\def\H{\mathcal H}
\def\A{\mathcal A}
\def\C{\mathcal C}
\def\D{\mathcal D}
\def\E{\mathcal E}
\def\G{\mathcal G}
\def\B{\mathcal B}
\def\L{\mathcal L}
\def\U{\mathcal U}
\def\K{\mathcal K}
\def\El{\mathcal E{\textit{l}}}

\def\M{\mathcal M}
\def\O{\mathcal O}
\def\R{\mathcal R}
\def\S{\mathcal S}
\def\N{\mathcal N}
\def\P{\mathcal P}

\newcommand{\undertilde}[1]{\underset{\sim}{#1}}
\newcommand{\abs}[1]{{\lvert#1\rvert}}
\newcommand{\mC}[1]{\mathfrak{C}(#1)}
\newcommand{\sigInf}[1]{\Sigma^{\infty}{#1}}
\newcommand{\x}[4]{\underset{#1, #2}{ \overset{#3, #4} \prod }}
\newcommand{\mA}[2]{\textit{Add}^n_{#1, #2}}
\newcommand{\mAK}[2]{\textit{Add}^k_{#1, #2}}
\newcommand{\mAL}[2]{\textit{Add}^l_{#1, #2}}
\newcommand{\Mdl}[2]{\L_\infty}
\newcommand{\inv}[1]{#1^{-1}}
\newcommand{\Lan}[2]{\mathbf{Lan}_{#1}(#2)}

\newcommand{\del}{\partial}
\newcommand{\sCatO}{\mathcal{S}Cat_\O}
\newcommand{\FCgop}{\mathbf{F}\mC{N(\gop)}}
\newcommand{\hProd}{{\overset{h} \oplus}}
\newcommand{\hProdn}{\underset{n}{\overset{h} \oplus}}
\newcommand{\hProdk}[1]{\underset{#1}{\overset{h} \oplus}}
\newcommand{\map}{\mathcal{M}\textit{ap}}
\newcommand{\SMGS}[2]{\map_{\gSC}({#1},{ #2})}
\newcommand{\MGS}[2]{\underline{\map}_{\gSC}({#1},{ #2})}
\newcommand{\MGSBox}[2]{\underline{\map}^{\Box}_{\gSC}({#1},{ #2})}
\newcommand{\Aqcat}[1]{\underline{#1}^\oplus}
\newcommand{\Cat}{\mathbf{Cat}}
\newcommand{\Sp}{\mathbf{Sp}}
\newcommand{\SpStb}{\mathbf{Sp}^{\textit{stable}}}
\newcommand{\SpStr}{\mathbf{Sp}^{\textit{strict}}}
\newcommand{\Sspec}{\mathbb{S}}
\newcommand{\ud}[1]{\underline{#1}}
\newcommand{\inrt}{\mathbf{Inrt}}
\newcommand{\act}{\mathbf{Act}}
\newcommand{\StrSMHom}[2]{[#1,#2]_\otimes^{\textit{str}}}
\newcommand{\Sh}[1]{{#1}^\sharp}
\newcommand{\Fl}[1]{{#1}^\flat}
\newcommand{\Nt}[1]{{#1}^\natural}
\newcommand{\Flmap}[3]{\Fl{\left[#1, #2\right]}_{#3}}
\newcommand{\Shmap}[3]{\Sh{\left[#1, #2 \right]}_{#3}}
\newcommand{\mRN}[1]{\int^{n^+ \in \gop}{#1}^+}
\newcommand{\mRNGen}[2]{\int_{+}^{\MakeLowercase{{#2}} \in {#2}}{#1(\MakeLowercase{{#2}})}}
\newcommand{\mRNL}[1]{{\mathfrak{F}}^+_{#1}(\gop)}
\newcommand{\mRNGenL}[2]{{\mathfrak{F}}^+_{#1}(#2)}
\newcommand{\mapG}[2]{[#1, #2]_{\gop}^+}
\newcommand{\mapFl}[2]{\Fl{[#1, #2]}}
\newcommand{\mapSh}[2]{\Sh{[#1, #2]}}
\newcommand{\mapMS}[2]{[#1, #2]^+}
\newcommand{\ExpG}[2]{{\left({#2}\right)}^{[#1]}}
\newcommand{\expG}[2]{{{#2}}^{[#1]}}
\newcommand{\mapGen}[3]{[#1, #2]_{#3}}
\newcommand{\mapGenM}[3]{[#1, #2]^+_{#3}}
\newcommand{\rNGen}[2]{\int^{\MakeLowercase{{#2}} \in {#2}}{#1}(\MakeLowercase{{#2}})}
\newcommand{\fRNGen}[2]{\fint^{\MakeLowercase{{#2}} \in {#2}}{#1}}
\newcommand{\rrNGen}[2]{\left(\rNGen{#1}{#2}\right)_\bullet}
\newcommand{\rNGenL}[2]{\mathfrak{L}_{#2}(#1)}
\newcommand{\rN}[1]{\int^{n^+ \in \gop} {#1}}
\newcommand{\rNGenR}[1]{\textit{h}^\ast}
\newcommand{\mRNGenR}[1]{\textit{h}_{+}^\ast}
\newcommand{\cn}[1]{\mathbf{id}(d)_{{#1}}}
\newcommand{\cnGen}[2]{\mathbf{id}(#2)_{{#1}}}
\newcommand{\piSSGen}[1]{\left( \sSetsMGen{#1} \right)^{\pi_0}}
\newcommand{\colim}[1]{{\varinjlim}^{#1}}
\newcommand{\hColim}[1]{{\varinjlim}^{\textit{h}}_{#1}}

\begin{abstract}
	The main objective of this paper is to show that the homotopy colimit of a diagram of quasi-categories and indexed by a small category is a localization of Lurie's higher Grothendieck construction of the diagram. We thereby generalize Thomason's classical result which states that the homotopy colimit of a diagram of categories has the homotopy type of (the classifying space of) the Grothendieck construction of the diagram of categories.

\end{abstract}

\maketitle


\section[Introduction]{Introduction}
\label{Introduction}

 In the paper \cite{T79} it was shown that the homotopy type of the homotopy colimit of a diagram of categories is encoded in the \emph{Grothendieck} construction of the diagram. We first recall this result:
 \begin{thm} \cite[Thm. 1.2]{T79}
 	\label{Ho-Colim-Thm-Cat}
 	Let $F:D \to \Cat$ be a functor taking values in categories. There is a natural homotopy equivalence of spaces:
 	\begin{equation*}
  \lvert hocolim{ \ N F} \rvert  \overset{\sim} \to  \B\left(\rNGen{F}{D}\right),
 	\end{equation*}
 	between (geometric realization of) the homotopy colimit of $N F$ and the classifying space of the Grothendieck construction of $F$.
 \end{thm}
  This result is commonly referred to as the \emph{homotopy colimit theorem}.
In this paper we extend this well known result by formulating and proving a homotopy colimit theorem for diagrams indexed by a small category and taking values in quasi-categories. We recall that quasi-categories are a simplicial sets based model for $(\infty, 1)$-categories wherein  a composite of a pair of composable arrows is only required to be unique upto a contractible space of choices. 
Our homotopy colimit theorem has two equivalent versions in two different settings. The first version is stated in the setting of marked simplicial sets. We recall that a \emph{marked} simplicial set is a pair $(X, \E)$ consisting of a simplicial set $X$ and a set $\E$ of edges ($1$-simplices) of $X$ which contains all degenerate edges. Each simplicial set can be \emph{naturally marked} by considering all \emph{equivalences}, see definition \ref{eq-sSet}, of the simplicial set as the marked edges. Now we state the first version of our homotopy colimit theorem:
\begin{thm}
	\label{Ho-Colim-Thm-QCat-M}
	Let $F:D \to \QCat$ be a functor taking values in quasi-categories. There is a natural homotopy equivalence of marked quasi-categories
	\begin{equation*}
		\emph{hocolim}{ \ F_+} \overset{\sim} \to  \mRNGen{F_+}{D} ,
	\end{equation*}
	between Lurie's marked (higher) Grothendieck construction of $F_+$ which is a naturally marked version of $F$ obtained by marking all equivalences in the quasi-category $F(d)$, for each $d \in D$, and the homotopy colimit of $F_+$.
\end{thm}
 The category of marked simplicial sets and simplicial maps preserving the marked edges is denoted $\sSetsM$.
  The category of  simplicial sets is related to that of marked simplicial sets by an adjunction:
\begin{equation*}
	L:\sSetsM \rightleftarrows \sSets:E
\end{equation*}
The right adjoint $E$ is the functor which \emph{marks} those edges in a simplicial set which are equivalences. The left adjoint $L$ is a localization functor which \emph{inverts} all marked edges of a marked simplicial set.
The second version of our homotopy colimit theorem is stated in the setting of (unmarked) simplicial sets as follows:
\begin{thm}
	\label{Ho-Colim-Thm-QCat}
	Let $F:D \to \QCat$ be a functor taking values in quasi-categories. There is a natural homotopy equivalence of (unmarked) quasi-categories:
	\begin{equation*}
		\emph{hocolim}{ \ F} \overset{\sim} \to L\left( \mRNGen{F_+}{D}  \right),
	\end{equation*}
	between the homotopy colimit of $F$ and a localization of Lurie's marked (higher) Grothendieck construction of $F_+$ which is a naturally marked version of $F$ obtained by marking all equivalences in the quasi-category $F(d)$, for each $d \in D$.
\end{thm}
 
 Clearly our theorem applies to a much larger class of diagrams but the two theorems are based in different settings.
  Unlike Thomason's homotopy colimit theorem, which merely determines the homotopy type of the homotopy colimit of a diagram of categories, our homotopy colimit theorem determines the homotopy colimit of a diagram of quasi-categories, up to equivalence of quasi-categories. Thus the two theorems may not agree even for a diagram of categories.
 To illustrate the distinction let us consider a functor $F:D \to \Cat$.
 Thomason's theorem says that $hocolim \ N F$ has the homotopy type of the classifying space of the Grothendieck construction of $F$. On the other hand, our theorem says that $hocolim \ NF$ is categorically equivalent to (the nerve of) the localization of the Grothendieck construction of $F$ obtained by inverting only the \emph{horizontal} arrows in the Grothendieck construction of $F$.
 
  The category of marked simplicial sets $\sSetsM$ inherits a model category structure from the Joyal model category of simplicial sets $\sSetsQ$, see \cite{AJ1}, \cite{AJ2}, which we denote by $\sSetsMQ$. A marked simplicial set is fibrant in this model structure if and only if it is a naturally marked quasi-category \emph{i.e.} the simplicial set is a quasi-category and its set of marked edges consist of equivalences. The adjunction $(L, E)$ described above is a Quillen equivalence between $\sSetsMQ$ and $\sSetsQ$. We will further review marked simplicial sets in Appendix \ref{mar-sSets}.
 
   In \cite[Sec. 3.2.5]{JL}, the author constructs a \emph{relative nerve} of a diagram of simplicial sets indexed by a (small) category $D$. Throughout this paper, except in Appendix \ref{comp-rel-ner}, we will refer to the relative nerve as \emph{Lurie's higher Grothendieck construction} or simply as \emph{Lurie's Grothendieck construction} when the context is clear. This construction defines a right Quillen functor of a Quillen equivalence which we denote, using a different notation than the one in \cite[Sec. 3.2.5.]{JL}, as follows:
\[
\rNGen{-}{D}:[D, \sSets] \to \ovCatGen{\sSets}{N(D)}.
\]
  A version of the aforementioned Quillen equivalence adapted to diagrams indexed by a (small) category $D$ and taking values in marked simplicial sets  is implied by \cite[Thm. 3.2.5.18]{JL}. More precisely, a version of Lurie's Grothendieck construction functor for \emph{marked} simplicial sets:
  \[
 \mRNGen{-}{D}:[D, \sSetsMQ] \to \ovCatGen{\sSetsM}{N(D)}.
  \]
   is a right Quillen functor of a Quillen equivalence between the projective model category of diagrams taking values in the model category of (marked) quasi-categories $\sSetsMQ$ and the coCartesian model category whose fibrant objects are coCartesian fibrations over the nerve of $D$. A second prominent objective of this paper is to establish another Quillen equivalence between the two aforementioned model categories wherein the left Quillen functor goes in the opposite direction, namely from $[D, \sSetsMQ]$ to $\ovCatGen{\sSetsM}{N(D)}$. More precisely we want to show that a homotopy colimit functor is a left Quillen functor of a Quillen equivalence between the two aforementioned model categories.
   A standard construction of a homotopy colimit of a diagram of spaces $F:D \to \sSets$ is obtained by taking the \emph{diagonal} of the \emph{bar construction} on $F$. This construction defines a functor
   \begin{equation*}
   	\gCL{} : [D, \sSets] \to \sSetsGen{N(D)}.
   	\end{equation*}
   An $n$-simplex in (the total space of) $\gCL{}(F)$ is a pair $(\sigma, x)$, where $\sigma \in N(D)_n$ and $x \in F(\sigma(0))_n$. The functor $\gCL{}$ has a natural extension to diagrams taking values in marked simplicial sets:
   \begin{equation*}
   	\gCLM{} : [D, \sSetsM] \to \sSetsMGen{C}.
   \end{equation*}
    The underlying simplicial set of $\gCLM{}(F)$ is the same as that of $\gCLM{}(U(F))$, where $U$ is the obvious forgetful functor, and an edge $(\sigma, x)$ is marked in $\gCLM{}(F)$ if $x$ is marked in $F(\sigma(0))$. We show that the following composite functor:
    \begin{equation*}
    [D, \sSetsMQ] \overset{\gCLM{}}\to \sSetsMGen{C} \overset{u^\ast}\to \sSetsMQ
    \end{equation*}
     is a homotopy colimit functor.
    Another prominent result of this paper establishes a second Quillen equivalence between the aforementioned model categories whose \emph{left} Quillen functor is $\gCLM{D}$. We provide a direct and self contained proof of this result. Lurie's proof uses a lot of machinery and is different from ours because it does not exploit the relation between his (marked) Grothendieck construction functor and $\gCLM{D}$.
      In fact Lurie's Quillen equivalence follows from our result. We closely follow the approach presented in \cite{HM} to prove an analogous result in the context of the covariant model category structure and the projective model category structure with respect to the Kan-Quillen model category of simplicial sets. However we establish our results by proving a dual statement to the one proved in \cite{HM} and thereby imply a new proof of \cite[Thm. C]{HM}.
   
   A third objective of this paper is to critique Lurie's Grothendieck construction. We provide a new description of Lurie's Grothendieck construction of a simplicial diagram: From a simplicial diagram $F:D \to \sSets$ we extract a simplicial space and refer to the zeroth row of this simplicial space, viewed as a bisimplicial set, as Lurie's Grothendieck construction of $F$.
 We justify our terminology in Appendix \ref{comp-rel-ner} wherein we show that the aforementioned zeroth row is (naturally) ismorphic to Lurie's Grothendieck construction as described in \cite[Sec. 3.2.5.]{JL}.
Our point of view on the notion of a higher Grothendieck construction is that the higher simplicies should be obtained by a procedure which is an extension of the one used to obtain the objects and arrows of the classical Grothendieck construction, namely the vertices of $\rNGen{F}{C}$ are pairs $(c, \beta_0) \in Ob(D) \times F(D)_0$ and edges are triples $(f,\beta_0, \beta_1) \in Mor(D) \times F(s(f))_0 \times F(t(f))_1$ such that $F(f)(\beta_0) = d_1(\beta_1)$. 
 The idea of higher simplices is expressed by the following picture which depicts a $3$-simplex in $\rNGen{F}{C}$ over a $3$-simplex $(f_1, f_2, f_3) \in N(D)_3$:
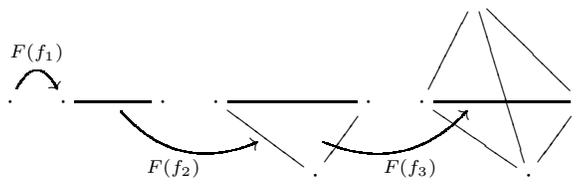
\begin{equation*}
\xymatrix@R=2mm@C=4mm{
	&&&&&& & && \cdot \ar@{-}[rrddd] \ar@{-}[rddddd] \ar@{-}[lddd] \\
&&&&&& & \\
&&&&&& & \\
\cdot \ar@/^1pc/[r]^{F(f_1)} & \cdot \ar@{-}[rr] &\ar@/_1pc/[rrrd]_{F(f_2)}& \cdot & \cdot \ar@{-}[rrr] \ar@{-}[rrdd] &&& \cdot \ar@{-}[ldd] & \cdot \ar@{-}[rrr] \ar@{-}[rrdd] &&&  \cdot \ar@{-}[ldd] \\
&&&&&&\ar@/_1pc/[rrru]_{F(f_3)} & &&&& \\
&&&&&& \cdot &&&& \cdot
}
\end{equation*}

	In section \ref{inf-cat-Gr-const} we construct \emph{ higher mapping path spaces} whose vertices are the aforementioned higher simplices and which capture the higher coherence data inherent in a diagram of quasi-categories.
	We obtain an elegant description of these by identifying them with mapping spaces in a simplicial functor category. These mapping path spaces play a vital role in our description of Lurie's Grothendieck construction.
	
	In appendix \ref{Classic-GC} we review the classical Grothendieck construction of a diagram of categories and show that it is isomorphic to the horizontal structure of a double category determined by the diagram of categories.
	Our construction of Lurie's Grothendieck construction is carried out in Section \ref{inf-cat-Gr-const}. An isomorphism between our definition of Lurie's higher Grothendieck construction with the one given in \cite[Defn. 3.2.5.2]{JL} is established in appendix \ref{comp-rel-ner}. A rectification theorem for coCartesian fibrations is proved in section \ref{coCart-mdl-str}.
	 The functor $\gCL{D}$ is usually NOT a homotopy colimit functor if the underlying category of simplicial sets is endowed with the Joyal model category structure. The main objective of section \ref{hoColim-QCat} is to suitably modify the marked version of this functor, namely $\gCLM{D}$, to produce a homotopy colimit functor for diagrams taking value in $\sSetsQ$. In the same section we also prove the main result, namely Theorem \ref{Ho-Colim-Thm-QCat}.

   \begin{ack}
   	The author would like to thank André Joyal for having useful discussions on the subject and for sharing his views on the notion of a higher Grothendieck construction.
   	\end{ack}

 \section{Lurie's Grothendieck construction revisited}
\label{inf-cat-Gr-const}
In this section we revisit Lurie's notion of higher Grothendieck construction for simplicial diagrams introduced in \cite[Sec. 3.2.5]{JL} where it is called the relative nerve. 
As mentioned earlier, our point of view on Lurie's Grothendieck construction is that it's merely an extension of the classical Grothendieck construction to functors taking values in simplicial sets.
Based on this point of view, in this section we provide another equivalent description of Lurie's Grothendieck construction of a simplicial diagram.
Later in this section we define a functor
 \[
  \fRNGen{-}{C} :[D, \sSets] \to \ovCatGen{\sSSets}{\Delta[0] \Box N(D)},
 \]
 where $\sSSets$ is the category of bisimplicial sets $[\Delta^{op} \times \Delta^{op}, \Sets]$ and the base bisimplicial set is defined as follows:
 \[
  (\Delta[0]\Box N(D))_{m,n} :=  \Delta[0]_m \times N(D)_n.
  \]
   We refer to the following composite as \emph{Lurie's (higher) Grothendieck construction} functor:
 \begin{equation*}
	 [D, \sSets] \overset{ \fRNGen{-}{C} } \to \ovCatGen{\sSSets}{N(D)\Box \Delta[0]} \overset{i_1^\ast} \to \ovCatGen{\sSets}{N(D)},
	 \end{equation*}
 where $i_1^\ast$ is the functor which restricts a bisimplicial set to its zeroth row.
We justify our terminology in Appendix \ref{comp-rel-ner} wherein we show that the aforementioned zeroth row is (naturally) ismorphic to Lurie's Grothendieck construction as described in \cite[Sec. 3.2.5.]{JL}.

The classical Grothendieck construction functor assigns to a functor $F:D \to \Cat$ an opfibration $p:\rNGen{F}{D} \to D$. The morphisms of the (total category) $\rNGen{F}{D}$ can be identified with the object set of the following category:
\begin{equation*}
Mor\left( \rNGen{F}{D} \right) := \underset{f \in Mor(D)} \bigsqcup \limproj{\A^1_{f}(F)}
\end{equation*}
where $\A^1_{f}(F)$ is the following zig zag:
\[
F(s(f)) \overset{F(f)} \to F(t(f)) \overset{d_1} \leftarrow F(t(f))^I
\]
The limit of the above diagram is referred to as the \emph{mapping path category} of the functor $F(f)$. The mapping path category is equipped with the following two functors which determine a category object in $\Cat$, also known as a \emph{double category}:
\begin{equation*}
\underset{f \in Mor(D)}\bigsqcup \limproj{\A^1_{f}(F)} \overset{s}{\underset{t} \rightrightarrows} \underset{c \in Ob(D)} \bigsqcup F(D)
\end{equation*}
The classical Grothendieck construction of $F$ is obtained by passing to the \emph{horizontal structure} of the double category which is obtained by applying the object function functor to the above diagram in $\Cat$ and passing to the category $\Sets$. In other words, the underlying graph of the category $\rNGen{F}{C}$ is the following diagram:
\begin{equation*}
\underset{f \in Mor(D)}\bigsqcup Ob(\limproj{\A^1_{f}(F)}) \overset{Ob(s)}{\underset{Ob(t)} \rightrightarrows} \underset{c \in Ob(D)} \bigsqcup Ob(F(D))
\end{equation*}

There is an analogous notion of a \emph{mapping path space} associated to a simplicial diagram. However, the mapping path space does not capture higher coherence data inherent in a finite sequence of simplicial maps.
 In order to capture this higher coherence data, for a pair $(f_1, f_2) \in N(D)$, we
 have to enhance the above zig zag to the following:
\[
F(c_0) \overset{F(f_1)} \to F(c_1) \overset{d_1} \leftarrow F(c_1)^{\Delta[1]} \overset{F(f_2)} \to F(c_2) \overset{d_2} \leftarrow F(c_2)^{\Delta[2]}
\]
 which we denote by $\A^2_{f_1,f_2}(F)$. We refer to the limit $\limproj{\A^2_{f_1,f_2}(F)}$ as the second mapping path space of $(F((f_1,f_2)))$ and denote it by $\P^2_{F((f_1, f_1))}$. The above discussion can be extended to a sequence $\sigma=(f_1, \dots, f_n)$ of $n$-composable arrows of $D$ which allows us to define the $n$th mapping path space of $F(\sigma)$ which we denote by $\P^n_{F(\sigma)}$. We show that the simplicial sets which are components of the following disjoint union:
 \[
\P^\bullet_F := \left\{ \underset{\sigma \in N(D)_n} \bigsqcup \P^n_{F(\sigma)}, n \in \Nat \right\}
 \]
  glue together into a simplicial space. Thus analogous to the classical Grothendieck construction which is obtained by passing to the horizontal structure (level of objects) of a category object in $\Cat$, Lurie's higher Grothendieck construction of a simplicial diagram is obtained by passing to the zeroth row (level of vertices) of a simplicial object in $\sSets$.

 Our new description facilitates a comparison between Lurie's Grothendieck construction of a diagram of simplicial sets with the standard homotopy colimit of the diagram (in the model category of spaces). In fact it allows us to construct a natural transformation between the two. Our description also facilitates in explicitly describing the structure maps of Lurie's Grothendieck construction which was left out in the original exposition. We provide a detailed account of these maps and also verify the simplicial identities. In order to describe our interpretation of Lurie's Grothendieck construction of a simplicial diagram, we first have to define a notion of a \emph{higher} mapping path space:
\begin{nota}
	We regard each object $[n]$ of the category $\Delta$ as a category depicted by the following linear graph:
	\[
	0 \to 1 \to \cdots \to n.
	\]
	We denote the unique arrow from $i$ to $i+1$ in the category $[n]$ by $(i, i+1)$, for $0 \le i \le n-1$.
	\end{nota}
The following two definitions will be used to define the desired notion of a higher mapping path space:
 \begin{df}
 	\label{nth-can-func-Path-Fib}
 For each $[n] \in \Delta$, there is a canonical functor $c(n):[n] \to \sSets$ which can be described by the following sequence:
 \begin{equation}
 \Delta[0] \overset{c_n((0,1))}  \to \Delta[1] \overset{c_n((1,2))}\to \cdots \overset{c_n((n-1,n))} \to \Delta[n],
 \end{equation}
 where the simplicial map $c_n((i-1,i))$ is \emph{face map} $d_i:\Delta[i-1] \to \Delta[i]$ and  $[n]$ is considered to be a category as described above. We refer to this canonical functor as the \emph{$n$th sequence of representables}.
 \end{df}
We recall from the definition of the nerve functor that an $n$-simplex in $N(D)$ is a functor $[n] \to C$.

\begin{df}
	For each $n$-simplex $\sigma$ in $N(D)$, a functor $F:D \to \sSets$ determines a sequence of length $n$ by the following composite:
	\[
	[n] \overset{\sigma} \to C \overset{F} \to \sSets.
	\]
	We refer to this composite functor as the \emph{$n$th sequence of $F$ over $\sigma$} and denote it by $F(\sigma)$.
	\end{df}
 
 The functor category $[[n], \sSets]$ inherits a simplicial category structure from the (simplicial) category of simplicial sets $\sSets$ \cite[Thm. 11.7.3]{Hirchhorn}.
 
  Now we are ready to define a higher mapping path space:
 \begin{df}
 	\label{n-path-Sp}
 	The \emph{$n$th mapping path space} of $X:[n] \to \sSets$, denoted by $\P_X^n$, is defined as follows:
 \begin{equation}
 \label{can-iso-map-PS}
 \P_X^n := \MapC{c(n)}{X}{[[n], \sSets]}
 \end{equation}
 \end{df}
\begin{nota}
	An $n$-simplex $\sigma $ in $N(D)$ is a functor $\sigma:[n] \to \sSets$ which can be described by the following diagram in $D$:
	\[
	c_0 \overset{f_1} \to c_1 \overset{f_2} \to c_2 \overset{f_3} \to \cdots \overset{f_{n-1}} \to c_{n-1} \overset{f_n} \to c_n
	\]
	\end{nota}
 Let $F:D \to \sSets$ be a functor and $\sigma$ be an $n$-simplex in $N(D)$. For this situation we will now make some remarks on the above definition:
\begin{rem}
	The $n$th mapping path space of $F(\sigma)$ is
	the limit of the following zig zag in $\sSets$:
	\begin{equation*}
	\xymatrix{
		&&& && F(c_n)^{\Delta[n]} \ar[d]^{F(c_n)^{d_n}} \\
		&&&&F(c_{n-1})^{\Delta[n-1]} \ar[r]_{F(f_n)^{\Delta[n-1]} \ \ } \ar[d]_{F(c_{n-1})^{d_{n-1}}} & F(c_n)^{\Delta[n-1]}   \\
		&&&F(c_{n-2})^{\Delta[n-2]} \ar[r]_{F(f_{n-1})^{\Delta[n-2]} \ \ } \ar[d]_{F(c_{n-2})^{d_{n-2}}} &F(c_{n-1})^{\Delta[n-2]} \\
		&&&F(c_{n-2})^{\Delta[n-3]} \ar@{--}[l] \\
		& F(c_2)^{\Delta[1]} \ar[d]^{F(c_2)^{d_1}} \ar[r]^{ \ \ F(f_2)^{\Delta[1]}}   & F(c_3)^{\Delta[1]} \ar@{--}[u] & \\
		F(c_1) \ar[r]_{F(f_1)} & F(c_2)
	} 
	\end{equation*}
\end{rem}
\begin{rem}
	\label{coord-rep}
	An $m$-simplex $\beta$ in $\P_{F(\sigma)}^n$, which is a map
	\[
	\beta: \TensP{\Delta[m]}{c(n)}{\sSets} \to F(\sigma)
	\]
	 in the functor category $[[n], \sSets]$, can be described as an $n$-tuple:
	\[
	(\beta_0, \beta_1, \dots, \beta_n) \in F(c_0)_m \times F(c_1)_{m+1} \times \cdots \times F(c_n)_{m+n}
	\]
	such that $\left(F(f_i)^{\Delta[m] \times {\Delta[i-1]}} \right)(\beta_{i - 1}) = d_{i}(\beta_{i})$, for $1 \le i \le n$.
	
\end{rem}
\begin{rem}
	\label{sig-one}
	When $\sigma$ is a one simplex $f:D_0 \to c_1$, the simplicial-set $\P_{F(f)}^1$ is the following pullback:
	\begin{equation*}
	\xymatrix{
		\P_{F(f)}^1 \ar[r] \ar[d] & F(c_1)^{\Delta[1]} \ar[d]^{F(c_1)^{d_1}} \\
		F(c_0) \ar[r]_{F(f)} & F(c_1)
	}
	\end{equation*}
	In other words, the simplicial set $\P_{F(f)}^1 = \MapC{c_1}{F(\sigma)}{[[1], \sSets]}$ is the classical mapping path space of the simplicial map $F(f):F(c_0) \to F(c_1)$. It follows immediately that
	\[
	\P^1_{F(id_c)} \cong [\Delta[1]; F(D)].
	\]
	Further, if $\sigma = (id_c, id_c, \dots, id_c)$ then
	
	\[
	\P^n_{F(\sigma)} \cong [\Delta[n]; F(D)].
	\]
\end{rem}


\begin{df}
	\label{codegen-can-seq}
	For each $n \ge 0$ and $0 \le i \le n$ we define a functor $s_i(c(n)):[n+1] \to \sSets$ which is represented by the following sequence:
	\begin{equation}
	\Delta[0] \overset{d_1} \to \Delta[1] \overset{d_2} \to \cdots \overset{d_{i}} \to \Delta[i] = \Delta[i] \overset{d_{i+1}} \to \Delta[i+1] \overset{d_{i+2}} \to \cdots \overset{d_{n}}\to \Delta[n]
	\end{equation}
	We will refer to this functor as the $i$th codegeneracy of $c(n)$.
\end{df}

The $(i-1)$th codegeneracy of $c(n-1)$ is equipped with a natural transformation
\begin{equation}
\label{face-map}
 d^{c(n-1)}_{i-1}:s_{i-1}(c(n-1)) \Rightarrow c(n)
\end{equation}
 which is represented by the following diagram:
\begin{equation*}
\xymatrix{
	\Delta[0] \ar[r]^{d_1} \ar@{=}[d] & \Delta[1] \ar[r]^{d_2} & \cdots \ar[r]^{d_{i-1}} & \Delta[i-1] \ar[r]^{d_i} & \Delta[i] \ar[r]^{d_{i+1}} & \Delta[i+1] \ar[r]^{d_{i+2}}  &\cdots \ar[r]^{d_n}& \Delta[n] \\
\Delta[0] \ar[r]_{d_1} & \Delta[1] \ar@{=}[u] \ar[r]_{d_2} & \cdots \ar[r]_{d_{i-1}} & \Delta[i-1] \ar@{=}[u] \ar@{=}[r] & \Delta[i-1] \ar[u]^{d_i} \ar[r]_{d_i} & \Delta[i] \ar[u]^{d_{i}} \ar[r]_{d_{i+1}} & \cdots \ar[r]_{d_{n-1}} & \Delta[n-1] \ar[u]^{d_i}
}
\end{equation*}

\begin{rem}
	\label{forgetful-fun}
	For each functor $X:[n] \to \sSets$, there is a canonical isomorphism
	\begin{equation*}
	U(n, X):\MapC{s_{i-1}(c(n-1))}{X}{[[n], \sSets]} \to \MapC{c(n-1)}{d_i(X)}{[[n-1], \sSets]},
	\end{equation*}
	where $d_i(X)$ is the following composite:
	\begin{equation*}
	[n - 1] \overset{d^i} \to [n] \overset{X} \to \sSets,
	\end{equation*}
	for $0 \le i \le n$. The map $U(n, X)$ maps a vertex $F:s_{i-1}(c(n-1)) \Rightarrow X$ to the following composite:
	\begin{equation*}
	\xymatrix{
		&  & \ar@{=>}[dd]^F &\\
	[n-1] \ar[r]^{d^i} & [n] \ar@/^3pc/[rr]^{s_{i-1}c(n)} \ar@/_3pc/[rr]_X && \sSets \\
	& &
   }
	\end{equation*}
	In the notation of remark \ref{coord-rep}, a $k$-simplex of the mapping space $\MapC{s_{i-1}(c(n-1))}{X}{[[n], \sSets]}$ can be represented by a pair $(\sigma, \beta)$, where $\sigma \in N(D)_n$ and $\beta$ is an $(n+1)$-tuple
	\begin{multline*}
	\beta = (\beta_o,\beta_1 \dots, \beta_{i-1}, X(i-1,i)(\beta_{i-1}), \dots, \beta_{n-1}) \\ \in X(0)_k \times X(1)_{k+1} \times \cdots \times X(i-1)_{k+i} \times X(i)_{k+i} \times \cdots X(n)_{k+n-1}.
	\end{multline*}
	In this notation, the simplicial map $U(n, X)$ can be described as follows:
	\[
	U(n, X)_k((\beta_o, \beta_1, \dots, \beta_{i-1}, X(i-1,i)(\beta_{i-1}), \dots, \beta_{n-1}) ) = (\beta_o, \beta_1, \dots, \beta_{i-1},  \dots, \beta_{n-1}) 
	\]
	\end{rem}
 
Each $n$-Path space $\P^n_{F(\sigma)}$ is equipped with $n+1$ \emph{face operators} :
\[
\P^n_{F(\sigma)}(d^i):\P^n_{F(\sigma)} \to \P^{n -1}_{F(d^{i}(\sigma))},
\]
for $0 \le i \le n$. The (simplicial) map $\P^n_{F(\sigma)}(d^i)$ is defined to be the following composite:
\begin{multline*}
\MapC{c(n)}{X}{[[n], \sSets]} \overset{\MapC{d^{c(n)}_i}{F(\sigma)}{[[n], \sSets]}} \to \MapC{s_{i-1}(c(n-1))}{F(\sigma)}{[[n], \sSets]} \overset{U(n,F(\sigma))} \to \\ \MapC{c(n-1)}{F(d_i(\sigma))}{[[n-1], \sSets]}
\end{multline*}
Using the notation of remark \ref{coord-rep}, the $i$th face operator defined above can be described, in degree $m$, as follows for $0 \le i \le n$:
\[
\left(\P_{F(\sigma)}^n(d^i)\right)_m ((\beta_0, \beta_1, \dots, \beta_n)) = (\beta_0, \beta_1, \dots, \beta_{i-1}, d_{i}(\beta_{i + 1}),   \dots,  d_{i}(\beta_{n-1}) , d_i(\beta_n)).
\]
A simple exercise in diagram chasing tells us the following:
\begin{multline*}
(\beta_0, \beta_1, \dots, \beta_{i-1}, d_{i}(\beta_{i + 1}),   \dots,  d_{i}(\beta_{n-1}) , d_i(\beta_n)) \in \\
F(\sigma(0))_m \times F(\sigma(1))_{m+1}, \times \cdots \times F(\sigma(i-1))_{m+i-1} \times F(\sigma(i+1))_{m+i} \times \cdots \times F(\sigma(n))_{m+n-1},
\end{multline*}
such that $F(\sigma((j-1, j)))(\beta_{j-1}) = d_j(\beta_j)$ for $1 \le j < i$, $F(\sigma((k-1, k)))(d_i(\beta_{k-1})) = d_{k-1}(d_i(\beta_{k}))$ for $i + 1 < k \le n$ and $F(d_i\sigma(i-1, i))(\beta_{i-1}) = d_i(d_i(\beta_{i+1}))$.
\begin{rem}
	In this remark we explicitly describe two extreme cases of the face maps:
The zeroth face map is described as follows:
\[
\left(\P_{F(\sigma)}^n(d^0)\right)_m ((\beta_0, \beta_1, \dots, \beta_n)) = (d_0\beta_1, \dots, d_{0}(\beta_{i}),   \dots,  d_0(\beta_{n-1}) , d_0(\beta_n)).
\]
 The $n$th face map is described as follows:
 \[
 \left(\P_{F(\sigma)}^n(d^n)\right)_m ((\beta_0, \beta_1, \dots, \beta_n)) = (\beta_0, \beta_1, \dots, \beta_{i},   \dots,  \beta_{n-1}).
 \]
 \end{rem}

The $i$th codegeneracy of $c(n)$ is equipped with a natural transformation
\begin{equation}
s^{c(n)}_{i}:D(n+1) \Rightarrow s_{i}(c(n))
\end{equation}
which is represented by the following diagram:
\begin{equation*}
\xymatrix{
	\Delta[0] \ar[r]^{d_1} \ar@{=}[d] & \Delta[1] \ar[r]^{d_2} & \cdots \ar[r]^{d_{i}} & \Delta[i] \ar[r]^{d_{i+1}} & \Delta[i+1] \ar[r]^{d_{i+2}} \ar[d]^{s_i} & \Delta[i+2] \ar[r]^{d_{i+3}} \ar[d]^{s_{i}}  &\cdots \ar[r]^{d_{n+1}}& \Delta[n+1] \ar[d]^{s_i} \\
	\Delta[0] \ar[r]_{d_1} & \Delta[1] \ar@{=}[u] \ar[r]_{d_2} & \cdots \ar[r]_{d_{i}} & \Delta[i] \ar@{=}[u] \ar@{=}[r] & \Delta[i]  \ar[r]_{d_{i+1}} & \Delta[i+1]  \ar[r]_{d_{i+2}} & \cdots \ar[r]_{d_{n}} & \Delta[n] 
}
\end{equation*}

The above natural transformation allows us to define $n+1$ \emph{degeneracy maps}
\[
\P_{F(\sigma)}^{n}(s^i):\P^{n}_{F(\sigma)} \to \P_{F(s_i(\sigma))}^{n + 1},
\]
for $0 \le i \le n$. The map $\P_F^n(s^i)$ is defined to be the following (simplicial) map:
\begin{multline*}
\label{degen-map-HGC}
\MapC{c(n)}{F(\sigma)}{[[n], \sSets]} \overset{\inv{U(n+1, F(s_i(\sigma)))}} \to \MapC{s_i({c(n)})}{F(s_i(\sigma))}{[[n+1], \sSets]} \\ \overset{\MapC{s_i^{c(n)}}{F(s_i(\sigma))}{[[n+1], \sSets]}} \to
\MapC{c(n+1)}{F(s_i(\sigma))}{[[n+1], \sSets]}
\end{multline*}
Using the notation of remark \ref{coord-rep}, the $i$th degeneracy operator defined above can be described, in degree $m$, as follows:
\[
\left( \P_F^n(s^i) \right)_m((\beta_0, \beta_1, \dots, \beta_{n-1}, \beta_n)) := (\beta_0, \dots, \beta_i, s_i(\beta_i), s_i(\beta_{i+1}), \dots, s_i(\beta_{n-1}), s_i(\beta_n)).
\]
A simple exercise in diagram chasing tells us the following:
\begin{multline*}
(\beta_0, \beta_1, \dots, \beta_{i}, s_i(\beta_i), s_{i}(\beta_{i + 1}),   \dots,  s_i(\beta_n)) \in \\
F(\sigma(0))_m \times F(\sigma(1))_{m+1}, \times \cdots \times F(\sigma(i))_{m+i} \times F(\sigma(i))_{m+i+1} \times \cdots \times F(\sigma(n))_{m+n+1},
\end{multline*}
such that $F(\sigma((j-1, j)))(\beta_{j-1}) = d_j(\beta_j)$ for $1 \le j \le i$ and $F(\sigma((k-1, k)))(s_i(\beta_{k-1})) = d_{k-1}(s_i(\beta_{k}))$ for $i + 1 \le k \le n$.

The following diagram depicts the image under $\P^2_{F((f_1, f_2))}(s^1)$ of a $0$-simplex $(\beta_0, \beta_1, \beta_2)$ in $\P^3_{F((f_1, id, f_2))}$
\[
\xymatrix@C=3mm@R=2.5mm{
&&&&&&&&&& \cdot \ar@{-}[rddddd] \ar@{-}[rrddd] \ar@{-}[lddd] \\
&&&&&&&&&&&&& \beta_2 \ar@/_1.8pc/[llld] \ar@/^1pc/[lld] \\
&&&&&&&&\ar@{}[d]_{ \ \ \ }&&&&& \\
\beta_0 \cdot \ar@/^1pc/[r]^{F(f_1)}&\cdot \ar@{-}[rr]^{\beta_1}&\ar@/_1pc/[rrrrd]_{F(id)}& \cdot&&\cdot \ar@{-}[rrr] \ar@{-}[rrdd] & &&\cdot \ar@{=}[ldd]& \ar@{-}[rrr]  \ar@{-}[rrdd]  &&& \cdot \\
&&&&&&&\ar@/_0.75pc/[rrru]_{F(f_2)} &&&&&& \\
&&&&&&&\cdot&&&& \cdot \ar@{=}[ruu] \\
&&&&&&&&&&&&& \\
&&&&&&&&&&&&& \\
&&&&&&\ar@{}[uuuu]_{s^1(\beta_1) \ \ \ }&&&&&&\ar@{}[uuuu]^{ s^1(\beta_2)}& \\
}
\]

 The classical Grothendieck construction defines a functor
\[
\rNGen{-}{C}:[C; \Cat] \to \ovCatGen{\Cat}{C} 
\]
We want to construct an extension of the above functor along the Nerve functor $[C; N]:[C; \Cat] \to [C; \sSets]$ which we call  \emph{Lurie's (higher) Grothendieck construction functor}. More precisely, we construct a functor
\[
\rNGen{-}{C}:[C; \sSets] \to \ovCatGen{\sSets}{N(D)} 
\]
such that the following diagram commutes up to isomorphism:
\begin{equation*}
\xymatrix@R=3mm@C=3mm{
[C; \Cat] \ar[r] \ar[dd] \ar@{=>}[rd]_\cong & \ovCatGen{\Cat}{C} \ar[r]^N & \ovCatGen{\sSets}{N(D)} \\
&&&\\
[C; \sSets] \ar@/_0.5pc/[rruu]_{\rNGen{-}{C}} &
}
\end{equation*}

\begin{df}
	The space of all $n$-Path spaces of $F:D \to \sSets$, denoted $\left( \fint^{c \in C} F \right)_n$, is defined to be the following simplicial set:
	\[
	\left( \fint^{c \in C} F \right)_n := \underset{\sigma \in N(D)_n} \bigsqcup   P^n_{F(\sigma)}.
	\]
	\end{df}

\begin{prop}
The spaces of all Path spaces of $F$ glue together into a simplicial space i.e. a functor $\fint^{c \in C} F:\Delta^{op} \to \sSets$ whose degree $n$ simplicial-set is defined as follows:
\[
\left( \fint^{c \in C} F \right)([n]):= \left( \fRNGen{F}{C} \right)_n
\]
 \end{prop}
\begin{proof}
	We begin by defining the degeneracy and face operators. The $i$th face operator $d_i:\left( \fRNGen{F}{C} \right)_n \to \left( \fRNGen{F}{C} \right)_{n-1}$ is defined as follows:
	\begin{equation*}
	d_i := \underset{\sigma \in N(D)_n} \bigsqcup P^n_{F(\sigma)}(d^i):\underset{\sigma \in N(D)_n} \bigsqcup   P^n_{F(\sigma)} \to \underset{\zeta \in N(D)_{n-1}} \bigsqcup   P^{n - 1}_{F(\zeta)}
	\end{equation*}
	The $i$th degeneracy operator $s_i:\left( \fRNGen{F}{C} \right)_{n} \to \left( \fRNGen{F}{C} \right)_{n+1}$ is defined as follows:
	\begin{equation*}
	s_i := \underset{\sigma \in N(D)_{n}} \bigsqcup P^{n}_{F(\sigma)}(s^i):\underset{\sigma \in N(D)_{n}} \bigsqcup   P^{n}_{F(\sigma)} \to \underset{\zeta \in N(D)_{n+1}} \bigsqcup   P^{n+1}_{F(\zeta)}
	\end{equation*}
	Since the face and degeneracy maps defined above are simplicial maps, therefore it is sufficient to verify the simplicial identities \cite[Pg. 25 (***)]{GZ} on an arbitrary $k$-simplex.
	We will verify the simplicial identities using the notation of remark \ref{coord-rep}. A $k$-simplex in $\left( \fRNGen{F}{C} \right)_n$ is a pair $(\sigma, \beta)$, where $\sigma \in N(D)_n$ and $\beta$ is an $(n + 1)$-tuple:
	\[
	\beta = (\beta_o, \beta_1, \dots, \beta_n) \in F(\sigma)(1)_k \times F(\sigma)(2)_{k+1} \times \cdots \times F(\sigma)(n)_{k+n},
	\]
	such that $F(\sigma)((i-1, i))(\beta_{i-1}) = d_i(\beta_i)$, for $1 \le i \le n$.
	 We now verify the simplicial identity $d_id_j = d_{j-1}d_i \ \ (i < j)$:
	\begin{multline*}
	d_id_j((\beta_0, \beta_1, \dots, \beta_n)) = d_i((\beta_0, \dots, \beta_{j-1}, d_j(\beta_{j+1}), \dots, d_j(\beta_n))) = \\
	(\beta_0, \dots, \beta_{i-1}, d_i(\beta_{i+1}), \dots d_i(\beta_{j-1}), d_id_j(\beta_{j+1}), \dots, d_id_j(\beta_n))) = \\
	(\beta_0, \dots, \beta_{i-1}, d_i(\beta_{i+1}), \dots d_i(\beta_{j-1}), d_{j-1}d_i(\beta_{j+1}), \dots, d_{j-1}d_i(\beta_n))) = \\
	d_{j-1}((\beta_0, \dots, \beta_{i-1}, d_i(\beta_{i+1}), \dots d_i(\beta_{j-1}), d_i(\beta_{j+1}), \dots, d_i(\beta_n)))) = \\
	d_{j-1}d_i((\beta_0, \dots, \beta_{i-1}, \beta_{i+1}, \dots \beta_{j-1}, \beta_{j+1}, \dots, \beta_n))).
	\end{multline*}
	The next simplicial identity which we verify is $s_is_j = s_{j+1}s_i$  $(i \le j)$:
	\begin{multline*}
	s_is_j((\beta_0, \dots, \beta_n)) = (\beta_0, \dots, \beta_j, s_j\beta_j, \dots s_j\beta_n) = \\
	(\beta_0, \dots, \beta_i, s_i(\beta_i), \dots,s_i(\beta_j), s_i s_j(\beta_j), \dots, s_is_j(\beta_n)) = \\
	(\beta_0, \dots, \beta_i, s_i(\beta_i), \dots,s_i(\beta_j), s_{j+1} s_i(\beta_j), \dots, s_{j+1} s_i(\beta_n)) = \\
	s_{j+1}((\beta_0, \dots, \beta_i, s_i(\beta_i), \dots,s_i(\beta_j),  \dots,  s_i(\beta_n))) = 
	s_{j+1}s_i ((\beta_0, \dots, \beta_n)).
	\end{multline*}
	An argument very similar to the one used in the verification of the above two simplicial identities, verifies the remaining simplicial identities. To avoid repetition, we leave the remaining verification for the interested reader.

	\end{proof}
\begin{nota}
	\label{Box-prod-sSets}
	Each pair $(K, L)$ of simplicial sets defines a \emph{bisimplicial sets} \emph{i.e.} a functor
	\[
	K \Box L: \Delta^{op} \times \Delta^{op}  \to \Sets
	 \]
	as follows:
	\[
	K \Box L([m],[n]) := K_m \times L_n
	\]
	\end{nota}
\begin{nota}
Each simplicial space $Z:\Delta^{op} \to \sSets$ determines a bisimplicial set, also denoted by $Z$

\[
Z:\Delta^{op} \times \Delta^{op} \to \Sets
\]
by $Z([m], [n]) = (Z([m], -))_n$. Further we denote the following simplicial set by $i_1^\ast(Z)$:
\[
\Delta^{op} \overset{ \ -\times [0] }\to \Delta^{op} \times \Delta^{op} \overset{Z} \to \Sets
\]
\end{nota}
Now we can define the Lurie's Grothendieck construction of a simplicial diagram $X:D \to \sSets$ as follows:
\begin{df}
	\label{LGrC}
	The \emph{Lurie's Grothendieck construction} of $X$ is the following simplicial set:
\begin{equation}
\rNGen{X}{C} = i^\ast_1 \left( \fRNGen{X}{C} \right)
\end{equation}
\end{df}
\begin{rem}
	\label{rep-Gr-cons}
	The set of $n$-simplices of $\rNGen{X}{C}$ can be represented as follows:
	\[
	\left( \rNGen{X}{C} \right)_n = \underset{\sigma \in N(D)_n} \bigsqcup  \left(\P^n_{X(\sigma)}\right)_0.
	\]
	An $n$-simplex of $\rNGen{X}{C}$ can be described as a pair $(\sigma, \beta)$, where $\sigma \in N(D)_n$ and $\beta$ is a natural transformation $\beta:c(n) \Rightarrow X(\sigma)$. Equivalently, the natural transformation $\beta$ can be fully described by an $(n+1)$-tuple
	\[
	(\beta_0, \beta_1, \dots, \beta_n) \in \left(X(\sigma)(0)\right)_0 \times \left(X(\sigma)(1)\right)_1 \times \cdots \times \left(X(\sigma)(n)\right)_n,
	\]
	such that $X(\sigma)((i, i+1))(\beta_i) = d_i(\beta_{i+1})$, for $0 \le i < n$.
	\end{rem}
\begin{rem}
	\label{proj-over-Ner}
	The simplicial space $ \fRNGen{X}{D} $ is equipped with a map of simplicial spaces:
	\[
	p^X_\bullet: \fRNGen{X}{D}  \to \Delta[0] \Box N(D).
	\]
\end{rem}

\begin{rem}
	\label{simp-sp-coTensor}
	Unwinding definitions given above one can deduce that the $n$th row of 
	the simplicial space $ \fRNGen{X}{D} $ maybe identified with the cotensor of Lurie's Grothendieck construction of $X$ in the simplicial category $\ovCatGen{\sSets}{N(D)}$
	\[
	\left( \fRNGen{X}{D} \right)(-, n) \cong \rNGen{X^{\Delta[n]}}{D} \cong \left[\Delta[n], \rNGen{X}{D}\right]_{D},
	\]
	where the mapping object $\left[\Delta[n], \rNGen{X}{D} \right]_{D}$ is the cotensor of $\rNGen{X}{D}$ with the representable simplicial set $\Delta[n]$ in $\ovCatGen{\sSets}{N(D)}$ and $X^{\Delta[n]}$ is the cotensor of $X$ with $\Delta[n]$ in the simplicial functor category $[D, \sSets]$.
\end{rem}

 \section{Rectification of coCartesian fibrations}
\label{coCart-mdl-str}
In this section we will prove a \emph{rectification theorem} for \emph{coCartesian} fibrations of simplicial sets over the nerve of a small category $D$. It was described in \cite[Ch. 3]{JL} that coCartesian fibrations over the nerve of a category are classified by a homotopy coherent diagram taking values in a (higher category of) quasi-categories. The rectification theorem which we  prove in this section says that these homotopy coherent diagrams can be replaced by honest functors. A similar result first appeared in \cite[Thm. 3.2.5.18]{JL}, wherein Lurie's higher Grothendieck construction functor for marked simplicial sets, as described in the introduction, is a right Quillen functor of a Quillen equivalence between the \emph{coCartesian} model category $\sSetsMGen{D}$ and the projective model category $[D, \sSetsMQ]$. The main objective of this section is to construct another Quillen equivalence between the same two model categories wherein the left Quillen functors goes in the opposite directions. We recall that the standard construction of homotopy colimit for diagrams of spaces, namely the diagonal of the bar construction, defines a functor
\begin{equation*}
	\gCL{D}: [D, \sSetsK] \to \sSetsGen{N(D)}.
\end{equation*}
An $n$-simplex in the total space of $\gCL{D}(F)$ is a pair $(\sigma, x)$, where $F:D \to \sSets$ is a functor, $\sigma$ is an $n$-simplex in $N(D)$ and $x \in F(\sigma(0))_n$.
The main result of \cite{HM} says that $\gCL{D}$ is a left Quillen functor of a Quillen equivalence between the projective model category structure on $[D, \sSetsK]$ and the covariant model category structure on $\sSetsGen{N(D)}$.
A marked version of the above functor
\begin{equation*}
	\gCLM{D}: [D, \sSetsMQ] \to \sSetsMGen{D}
\end{equation*}
is obtained by marking those edges of (the total space of) $\gCL{D}(F)$, where $F:D \to \sSetsM$ is a functor, which are determined by marked edges in the image of $F$. We first show in Proposition \ref{hocolim-func-mar-SS} that the composite functor
 \begin{equation*}
 [D, \sSetsMQ] \overset{\gCLM{D}} \to \sSetsMGen{D} \overset{u^\ast}\to \sSetsMQ
 \end{equation*}
  is a homotopy colimit functor in the sense of definition \ref{hocolim-func}. The second prominent result of this paper, namely Theorem \ref{main-res}, shows that the functor $\gCLM{D}$ is the left Quillen functor of a Quillen equivalence between the projective model category structure on $[D, \sSetsMQ]$ and the coCartesian model category structure on $\sSetsMGen{D}$ whose fibrant objects may be identified with coCartesian fibrations over $N(D)$. We prove our main result by the showing that the (right) derived functors of the following two functors:
  \begin{equation*}
  	[D, \sSetsMQ] \overset{\mRNGen{D}{D}} \to \sSetsMGen{D} \overset{\mRNGenR{D} \ } \to [D, \sSetsMQ]
  	\end{equation*}
  are mutual inverses on one another. We do so by constructing a natural transformation $\eta^+:id_{[D, \sSetsMQ]} \Rightarrow \mRNGenR{D}({\mRNGen{-}{D}})$ for each functor $F:D \to \sSetsMQ$, and a zig-zag of natural transformations
   \begin{equation*}
   \mRNGen{\mRNGenR{D}(-)}{D} \Leftarrow \gCLM{D}({\mRNGenR{D}}(-)) \Rightarrow id.
   \end{equation*}
  We show that for each fibrant $F:D \to \sSetsQ$,  $\eta^+_F$ is a natural weak-equivalence and for each fibrant $X \in \sSetsMGen{D}$, we have the following zig-zag of weak equivalences in $\sSetsMGen{D}$:
    \begin{equation*}
   	\mRNGen{\mRNGenR{D}(X)}{D} \leftarrow \gCLM{D}({\mRNGenR{D}}(X)) \rightarrow X.
   \end{equation*}
  A key ingredient in the proof of our main result is the construction of a natural transformation
\begin{equation*}
	\iota_D:\gCLM{arg1} \Rightarrow \mRNGen{-}{D}.
	\end{equation*}
 The construction of $\iota_D$ was facilitated by our description of the relative nerve functor in the previous sections. For each fibrant $F$ in the projective model category $[D, \sSetsMQ]$, we show that the map $\iota_D$ is a weak equivalence in the coCartesian model category $\sSetsMGen{D}$.

%
%
%
We begin with a review of coCartesin fibrations over the simplicial set $N(D)$. We will also review a model category structure on the category $\sSetsMGen{D}$ in which the fibrant objects are (essentially) coCartesian fibrations.
\begin{df}
	\label{p-CC-edge}
	Let $p:X \to S$ be an inner fibration of simplicial sets. Let $f:x \to y \in (X)_1$ be an edge in $X$. We say that $f$ is $p$-coCartesian if, for all $n \ge 2$ and every (outer) commutative diagram, there exists a (dotted) lifting arrow which makes the entire diagram commutative:
	\begin{equation}
	\xymatrix{
	\Delta^{\lbrace 0, 1 \rbrace} \ar@{_{(}->}[d] \ar[rd]^f \\
	\Lambda^0[n] \ar@{_{(}->}[d] \ar[r] & X \ar[d]^p \\
	\Delta[n] \ar[r] \ar@{-->}[ru] & S
   }
	\end{equation}
	\end{df}
\begin{rem}
	Let $M$ be a (ordinary) category equipped with a functor $p:M \to I$, then an arrow $f$ in $M$, which maps isomorphically to $I$, is coCartesian in the usual sense if and only if $f$  is $N(p)$-coCartesian in the sense of the above definition, where $N(p):N(M) \to \Delta[1]$ represents the nerve of $p$.
	\end{rem}
This definition leads us to the notion of a coCartesian fibration of simplicial sets:
\begin{df}
	\label{CC-fib}
	A map of simplicial sets $p:X \to S$ is called a \emph{coCartesian} fibration if it satisfies the following conditions:
	\begin{enumerate}
		\item $p$ is an inner fibration of simplicial sets.
		\item for each edge $p:x \to y$ of $S$ and each vertex $\ud{x}$ of $X$ with $p(\ud{x}) = x$, there exists a $p$-coCartesian edge $\ud{f}:\ud{x} \to \ud{y}$ with $p(\ud{f}) = f$.
		\end{enumerate}
	\end{df}
 A coCartesin fibration roughly means that it is up to weak-equivalence determined by a \emph{functor} from $S$ to a suitably defines $\infty$-category of $\infty$-categories. This idea is explored in detail in \cite[Ch. 3]{JL}.

\begin{nota}
	To each coCartesian fibration $p:X \to N(D)$ we can associate a marked simplicial set denoted $\Nt{X}$ which is composed of the pair $(X, \E)$, where $\E$ is the set of $p$-coCartesian edges of $X$
	\end{nota}
\begin{nota}
	\begin{sloppypar}
Let $(X, p), (Y, q)$ be two objects in $\sSetsMGen{D}$. We denote by  $[X, Y]_D^+$, the full (marked) simplicial subset of $[X, Y]^+$ spanned by maps in  $\sSetsMGen{D}(X, Y)$, namely spanned by maps in $[X, Y]^+$ which are compatible with the projections $p$ and $q$. We denote by  $\Flmap{X}{Y}{D}$, the full simplicial subset of $\Fl{[X, Y]}$ spanned by maps in  $\sSetsMGen{D}(X, Y)$. We denote by $\Shmap{X}{Y}{D} \subseteq \Shmap{X}{Y}{}$ the simplicial subsets spanned by maps in $\sSetsMGen{D}$.
\end{sloppypar}
\begin{df}
	\label{CC-Eq}
	A morphism $F:X \to Y$ in the category $\sSetsMGen{D}$ is called a \emph{coCartesian}-equivalence if for each coCartesian fibration $p:Z \to N(D)$, the induced simplicial map
	\[
	\Flmap{F}{\Nt{Z}}{D}:\Flmap{Y}{\Nt{Z}}{D} \to \Flmap{X}{\Nt{Z}}{D}
	\]
	is a categorical equivalence of simplicial-sets(quasi-categories).
	\end{df}
\begin{prop}
	\label{char-cc-eq}
	Let $u:X \to Y$ be a map in $\sSetsMGen{D}$, then the following are equivalent
	\begin{enumerate}
	 \item $u$ is a coCartesian equivalence.
	 \item For each functor $Z:D \to \sSetsM$, such that $Z(d)$ is a quasi-category whose marked edges are equivalences, the following (simplicial) map is a categorical equivalence:
	\[
	\Flmap{u}{ \mRNGen{Z}{D}}{D}:\Flmap{Y}{\mRNGen{Z}{D}}{D} \to \Flmap{X}{\rNGen{Z}{D}}{D}
	\]
	\item For each functor $Z:D \to \sSetsM$, such that $Z(d)$ is a quasi-category whose marked edges are equivalences, the following map is a bijection:
	\[
	\pi_0\Shmap{u}{ \mRNGen{Z}{D}}{D}:\pi_0\Shmap{Y}{\mRNGen{Z}{D}}{D} \to \pi_0\Shmap{X}{\rNGen{Z}{D}}{D}
	\]
	\end{enumerate}
\end{prop}
\begin{proof}
	$(1 \Rightarrow 2)$ Follows from the definition of coCartesian equivalence because $\mRNGen{Z}{D}$ is a coCartesian fibration under the given hypothesis.

	 Let us assume that $\Flmap{u}{ \mRNGen{Z}{D}}{D}$ is a categorical equivalence of quasi-categories for each functor $Z$ satisfying the given hypothesis.  This implies that $\Flmap{u}{ \Nt{T}}{D}$ is a categorical equivalence if and only if $\Flmap{u}{ \mRNGen{Z(T)}{D}}{D}$ is one.
		
		 $(2 \Rightarrow 3)$ We recall from \cite[Prop. 3.1.3.3]{JL} and \cite[Prop. 3.1.4.1]{JL} that, for any coCartesian fibration $\Nt{T} \in \sSetsMGen{D}$, the simplicial map $\Flmap{u}{ \Nt{T}}{D}$ is a categorical equivalence if and only if the map $\Shmap{u}{\Nt{T}}{D}$ is a homotopy equivalence of Kan complexes. This implies that $\pi_0\Shmap{u}{ \mRNGen{Z}{D}}{D}$ is a bijection.
		 
		  $(3 \Rightarrow 1)$ 
		  We recall from \cite[Cor. 3.1.4.4]{JL} that the coCartesian model category is a simplicial model category with simplicial function object given by the bifunctor $\Shmap{-}{-}{D}$.
		  This implies that $u$ is a coCartesian equivalence if and only if $\pi_0\Shmap{u}{\Nt{W}}{D}$ is a bijection for each fibrant object $W$ of the coCartesian model category. By  \cite[Prop. 3.1.4.1]{JL} we may replace $W$ by a coCartesian fibration $W \cong \Nt{T}$.
		  Further, it follows from \cite[Prop. 3.2.5.18(2)]{JL}  that for each cocartesian fibration $\Nt{T}$ there exists a functor $Z(T):D \to \sSetsM$, which satisfies the assumptions of the functor in the statement of the proposition, such that there is map
		  	\[
		  	F_T:\Nt{T} \to \mRNGen{Z(T)}{D}
		  	\]
		  	which is a coCartesian equivalence. Now it follows that $u$ is a coCartesian equivalence if and only if $\pi_0\Shmap{u}{\mRNGen{Z(T)}{D}}{D}$ is a bijection for each functor $Z$ satisfying the conditions mentioned in the statement of the proposition.
		 
\end{proof}

\end{nota}
Next we will recall a model category structure on the overcategory $\sSetsMGen{D}$ from \cite[Prop. 3.1.3.7.]{JL} in which fibrant objects are (essentially) coCartesian fibrations.
\begin{thm}
	\label{CC-Mdl-Str}
	There is a left-proper, combinatorial model category structure on the category $\sSetsMGen{D}$ in which a morphism is 
	\begin{enumerate}
		\item a cofibration if it is a monomorphism when regarded as a map of simplicial sets.
		\item a weak-equivalences if it is a coCartesian equivalence.
		\item a fibration if it has the right lifting property with respect to all maps which are simultaneously cofibrations and weak-equivalences. 
		\end{enumerate}
	\end{thm}

We have defined a function object for the category $\sSetsMGen{D}$ above.
The simplicial set $\Flmap{X}{Y}{D}$ has vertices, all maps from $X$ to $Y$ in $\sSetsMGen{D}$. An $n$-simplex in $\Flmap{X}{Y}{D}$ is a map $\Fl{\Delta[n]} \times X \to Y$ in $\sSetsMGen{D}$, where $\Fl{\Delta[n]} \times (X, p) = (\Fl{\Delta[n]} \times X, pp_2)$, where $p_2$ is the projection $\Fl{\Delta[n]} \times X \to X$. The enriched category $\sSetsMGen{D}$ admits tensor and cotensor products. The \emph{tensor product} of an object $X = (X, p)$ in $\sSetsMGen{D}$ with a simplicial set $A$ is the objects
\[
\Fl{A} \times X = (\Fl{A} \times X, pp_2).
\]
The \emph{cotensor product} of $X$ by $A$ is an object of $\sSetsMGen{D}$ denoted $\expG{A}{X}$. If $q:\expG{A}{X} \to \Sh{N(D)}$ is the structure map, then a simplex $x:\Fl{\Delta[n]} \to \expG{A}{X}$ over a simplex $y = qx:\Delta[n] \to \Sh{N(D)}$ is a map $\Fl{A} \times (\Fl{\Delta[n]}, y) \to (X, p)$. The object $(\expG{A}{X}, q)$ can be constructed by the following pullback square in $\sSetsM$:
\begin{equation*}
\xymatrix{
\expG{A}{X} \ar[r] \ar[d]_q & \mapMS{\Fl{A}}{X} \ar[d]^{\mapMS{\Fl{A}}{p}} \\
\Sh{N(D)} \ar[r] & \mapMS{\Fl{A}}{\Sh{N(D)}}
}
\end{equation*}
where the bottom map is the diagonal. There are canonical isomorphisms:
\begin{equation}
\Flmap{\Fl{A} \times X}{Y}{D} \cong \left[A, \Flmap{X}{Y}{D} \right] \cong \Flmap{X}{\expG{A}{Y}}{D}
\end{equation}
\begin{rem}
	\label{Simp-Mdl-Cat}
	The coCartesian model category structure on $\sSetsMGen{D}$ is a simplicial model category structure with the simplicial Hom functor:
	\[
	\Shmap{-}{-}{D}:\sSetsMGen{D}^{op} \times \sSetsMGen{D} \to \sSets.
	\]
	This is proved in \cite[Corollary 3.1.4.4.]{JL}. The coCartesian model category structure is a $\sSetsQ$-model category structure with the function object given by:
	\[
	\Flmap{-}{-}{D}:\sSetsMGen{D}^{op} \times \sSetsMGen{D} \to \sSets.
	\]
	This is remark \cite[3.1.4.5.]{JL}.
	\end{rem}
\begin{rem}
	\label{enrich-mar-sSets}
	The coCartesian model category is a $\sSetsMQ$-model category with the Hom functor:
	\[
	\mapGenM{-}{-}{D}:\sSetsMGen{D}^{op} \times \sSetsMGen{D} \to \sSetsM.
	\]
	This follows from \cite[Corollary 3.1.4.3]{JL} by taking $S = N(D)$ and $T = \Delta[0]$, where $S$ and $T$ are specified in the statement of the corollary.
	\end{rem}
 \begin{df}
 	\label{mar-rel-Ner}
 	Let $F:D \to \sSetsM$ be a functor. We can compose it with the forgetful functor $U$ to obtain a composite functor $F:D \overset{F} \to \sSetsM \overset{U} \to \sSets$. Lurie's \emph{marked} Grothendieck construction of $F$, denoted $\mRNGen{F}{D}$, is the marked simplicial set $\left(\mRNGen{F}{D}, \E \right)$,
 	where the set $\E$ consists of those edges $\ud{e} = (e, h)$ of $\mRNGen{F}{D}$, see remark \ref{Rel-Ner-edge}, which determines a marked edge of the marked simplicial set $F(d)$, where $e:c \to d$ is an arrow in $D$.
 	\end{df}
 The above construction of the marked Grothendieck construction determines a functor 
 \begin{equation}
 \label{mark-Gr-const}
 \mRNGen{-}{D}:[D, \sSetsM] \to \sSetsMGen{D}. 
 \end{equation}
 The functor $ \mRNGen{-}{D}$ has a left adjoint which we denote by $\mRNGenL{\bullet}{D}$, see \cite[Rem. 3.2.5.5]{JL}. This functor is defined on objects as follows:
 \begin{equation*}
 \mRNGenL{X}{D}(d) = X \underset{\Sh{N(D)}} \times \Sh{N(D/d)},
 \end{equation*}
 where $(X, p)$ is an object in $\sSetsMGen{D}$.
 \begin{nota}
 	We will sometimes denote $\mRNGenL{X}{D}$ by $\mRNGenL{\bullet}{D}(X)$.
 	\end{nota}
 \begin{rem}
 	\label{left-adj-pres-fib-obj}
 	If $(X, p)$ is a fibrant object in $\sSetsMGen{D}$ then $\mRNGenL{X}{D}$ is a fibrant object in the projective model category $[D, \sSetsMQ]$.
 	\end{rem}

 Next we will define a marked version of the functor $\rNGenR{D}$, denoted $\mRNGenR{D}$:
 \begin{equation*}
 \mRNGenR{D}(X)(d) := [\Sh{N(\ovCatGen{d}{D})}, X]^+_D
 \end{equation*}
 where $X$ is an object of $\sSetsMGen{D}$. This functor has a left adjoint which we denote by $\gCLM{D}$.
 \begin{nota}
 	We denote by $\cn{n}$ the $n$-simplex of $N(D)$ represented by the constant functor $\cn{n}:[n] \to D$ having value $d$.
 	An $n$-simplex $\sigma$ in the simplicial set $N(d/D)$ is a natural transformation $\sigma:\cn{n} \Rightarrow \beta$, where $\beta$ is an $n$-simplex in $N(D)_n$ viewed as a functor 
 	$\beta:[n] \to D$.
   For each functor $F:D \to \sSetsM$, an $n$-simplex $\sigma$ in $N(d/D)$ determines a natural transformation $F \circ \sigma$ which is represented by the following commutative diagram:
 	\begin{equation*}
 		\xymatrix@C=16mm{
 		F(d) \ar@{=}[r] \ar@{=}[d] & F(d) \ar@{=}[r] \ar[d] & \cdots \ar@{=}[r] & F(d) 	   \ar@{=}[r] \ar[d] &F(d) \ar[d]  \\
 		F(d) \ar[r]_{F(\sigma(0, 1)) \ \ } & F(\sigma(1)) \ar[r]_{F(\sigma(1, 2)) } & \cdots \ar[r]_{F(\sigma(n-1, n)) \ \ } & F(\sigma(n)) 	   \ar[r]_{F(\sigma(n,n+1)) \ \ \ } &F(\sigma(n+1)) 
        }
 		\end{equation*}
 	We denote this natural transformation by $c(F, \sigma)$.
 	\end{nota}
 
 \begin{df}
 	\label{unit-comp-Gr}
 	Let $F:D \to \sSetsM$ be a functor.
 	For each $d \in D$ we define a map of marked simplicial sets 
 	\[\eta^+_F(d):F(d) \to [N(\ovCatGen{d}{D}), \mRNGen{F}{D}]^+_D.
 	\]
 	Let $x \in F(d)_n$ be an $n$-simplex in $F(d)$. This $n$-simplex defines a canonical map $\eta^+_F(d)(x):N(\ovCatGen{d}{D}) \times \Delta[n] \to \mRNGen{F}{D}$ in $\sSetsMGen{D}$ whose value on $(\cn{n}, id_n) \in (N(\ovCatGen{d}{D}) \times \Delta[n])_n$ is the pair $(\cn{n}, \unit{}(x)) \in \left(\mRNGen{F}{D}\right)_n$, where $\unit{}(x):c(n) \Rightarrow F(\cn{n})$ is the natural transformation which is represented by the following diagram:
 	\begin{equation*}
 			\xymatrix{
 			\Delta[0] \ar[r]^{d_1} \ar[d] & \Delta[1] \ar[r]^{d_2} \ar[d] & \cdots \ar[r]^{d_{n-1}} & \Delta[n-1] \ar[r]^{d_n} \ar[d]_{} & \Delta[n] \ar[d]^{x} \\
 			F(d) \ar@{=}[r] &  F(d) \ar@{=}[r] & \cdots \ar@{=}[r]   &  F(d) \ar@{=}[r] &  F(d)
 		}
 		\end{equation*}

 	 We recall that a $k$-simplex in $\Delta[n]$ is a map $\alpha:[k] \to [n]$ in the category $\Delta$ and therefore it can be written as $\Delta[n](\alpha)(id_n)$.
 	For a $k$-simplex $(\sigma, \alpha)$ in $N(\ovCatGen{d}{D}) \times \Delta[n]$, where $\sigma = ((g, f_1, \dots, f_{k}) \in N(d/D)_k$ , we define the $k$-simplex $\eta^+_X(d)(x)((g, f_1, f_2, \dots, f_{k+1}), \alpha)$ to be the pair $(d_0(\sigma), c(F, \sigma) \cdot \unit{}(x))$, where the composite natural transformation $c(F, \sigma) \cdot \unit{}(x)$ is depicted by the following diagram:
 		\begin{equation*}
 		\xymatrix{
 			\Delta[0] \ar[r]^{d_1} \ar[d]_{} & \Delta[1] \ar[r]^{d_2} \ar[d]_{} & \cdots \ar[r]^{d_{k-1}} & \Delta[k-1] \ar[r]^{d_k} \ar[d]_{} & \Delta[k] \ar[d]^{F(d)(\alpha)(x)} \\
 			F(d) \ar@{=}[r] \ar@{=}[d]  &  F(d) \ar@{=}[r] \ar[d] & \cdots \ar@{=}[r]   &  F(d) \ar@{=}[r] \ar[d]  &  F(d) \ar[d]^{F(f_{k} \circ f_{k-1} \circ \cdots \circ f_1  \circ g)}  \\
 			F(d) \ar[r]_{F(g)} &  F(d_1) \ar[r]_{F(f_1)} & \cdots \ar[r]_{F(f_{k-1}) \ \ \ }   &  F(d_{k-1}) \ar[r]_{F(f_{k})} &  F(d_{k})
 		}
 	\end{equation*}
%
 	This defines the (simplicial) map $\eta^+_F(d)(x)$. These simplicial maps glue together into a natural transformation $\eta^+_F$.
 	\end{df}
%
 Now we define a map $\iota_d^+$ in $\sSetsMGen{D}$:
 \begin{equation}
 \label{main-gen-local-mar}
 \xymatrix{
 	\Fl{\Delta[0]} \ar[rr]^{id_d} \ar[rd]_d  && \Sh{N(\ovCatGen{d}{D})} \ar[ld] \\
 	&\Sh{N(D)}
 }
 \end{equation}
 \begin{lem}
 	\label{main-lemma-mar}
 	For each $d \in D$ the morphism $\iota^+_d$ defined in \eqref{main-gen-local-mar} is a coCartesian equivalence.
 \end{lem}
 \begin{proof}
 	We will show that for each functor $Z:D \to \sSetsM$ such that, for each $d \in D$, $Z(d)$ is a quasi-category whose marked edges are equivalences, we have the following bijection:
 	\[
 	\pi_0\Shmap{\iota^+_d}{\rNGen{Z}{D}}{D}:\pi_0\Shmap{N(\ovCatGen{d}{D})}{\rNGen{Z}{D}}{D}  \to \pi_0\Shmap{\Delta[0]} {\rNGen{Z}{D}}{D} \cong \pi_0(J(Z(d))),
 	\]
 	where $J(Z(d))$ is the largest Kan complex contained in $Z(d)$.
 	Let $z \in J(Z(d))_0$ be a vertex of $\Sh{J(Z(d))}$. We will construct a morphism $F_z:N(\ovCatGen{d}{D}) \to \rNGen{Z}{D}$ in the category $\sSetsMGen{D}$.  The vertex $z$ represents a natural transformation
 	\[
 	T_z:D(d, -) \Rightarrow Z
 	\]
 	such that $T_z(id_d) = z$. Since $N(\ovCatGen{d}{D}) \cong \rNGen{D(d, -)}{D}$ therefore we have a map
 	\[
 	F_z:N(\ovCatGen{d}{D}) \cong \rNGen{D(d, -)}{D} \overset{\rNGen{T_z}{D}} \to \rNGen{Z}{D}
 	\]
 	in $\sSetsMGen{D}$ such that $F_z(id_d) = z$. Thus we have shown that the map $\pi_0\Shmap{\iota^+_d} {\rNGen{Z}{D}}{D}$ is a surjection.
 	
 	Let $f:y \to z$ be an edge of $J(Z(d))$, then by the (enriched) Yoneda's lemma followed by an application of the Grothendieck construction functor, this edge uniquely determines a map
 	\[
 	T_f:N(\ovCatGen{d}{D}) \times \Delta[1] \to \rNGen{Z}{D}
 	\]
 	in $\sSetsMGen{D}$ such that $F_z((id_d, id_1)) = f$. Thus we have shown that the map $\pi_0\Shmap{\iota^+_d} {\rNGen{Z}{D}}{D}$ is also an injection.
 	
 \end{proof}

\begin{lem}
	\label{unit-Eq-mar}
	For any projectively fibrant functor $F:D \to \sSetsM$, the  map $\eta^+_F$ defined in \ref{unit-comp-Gr} is an objectwise categorical equivalence of marked simplicial sets.
\end{lem}
\begin{proof}
	Under the hypothesis of the lemma, it follows from \cite[Prop. 3.2.5.18(2)]{JL} and lemma \ref{isom-Rel-Ner} that $\mRNGen{F}{D}$ is a fibrant object in the coCartesian model category. Now lemma \ref{main-lemma-mar} and remark \ref{enrich-mar-sSets} gives us, for each $d \in D$, the following homotopy equivalence in $\sSetsMQ$:
	\[
	[\iota^+_d, \rNGen{F}{D}]^+_D:[\Sh{N(\ovCatGen{d}{D})}, \mRNGen{F}{D}]^+_D \to [\Fl{\Delta[0]}, \mRNGen{F}{D}]^+_D
	\]
	such that $c \circ [\iota_d, \rNGen{F}{D}]_D \circ \eta_F(d) = id_{F(d)}$, where 
	\begin{equation*}
		c:[\Fl{\Delta[0]}, \mRNGen{F}{D}]^+_D \cong F(d)
		\end{equation*}
	is the canonical isomorphism between the fiber of $p:\mRNGen{F}{D} \to N(D)$ over $d \in D$ and $F(d)$ \emph{i.e.} the value of the functor $F$ on $d$.
	Now the $2$ out of $3$ property of weak equivalences in a model category tells us that $\eta^+_F(d)$ is a homotopy equivalence for each $d \in D$ therefore $\eta^+_F$ is an objectwise homotopy equivalence in $[D,\sSetsMQ]$.
\end{proof}
 
 An immediate consequence of the definition of the right adjoint functor $\mRNGenR{D}$ is the following lemma:
 \begin{lem}
 	\label{Qu-Eq-L-R-mar}
 	The adjunction $(\gCLM{D}, \mRNGenR{D})$ is a Quillen adjunction between the projective model category structure on $[D, \sSetsMQ]$ and coCartesian model category $\sSetsMGen{D}$.
 \end{lem}
 \begin{proof}
 	The coCartesian model category is a $\sSetsMQ$-model category, see remark \ref{enrich-mar-sSets}. This implies that $\mRNGenR{D}$ maps (acyclic) fibrations in the coCartesian model category to (acyclic) projective fibrations in $[D, \sSetsM]$ which are objectwise (acyclic) fibrations of marked simplicial sets.
 \end{proof}
We have shown the existence of the left Quillen functor $\gCLM{D}$ above. Now we will provide a construction of this left adjoint. Let $F:D \to \sSetsM$ be a functor. An $n$-simplex in the (total space of) $\gCLM{D}(F)$ is a pair $(\sigma, x)$, where $\sigma = (d_0 \overset{f_1} \to d_1 \overset{f_1} \to d_2 \overset{f_2} \to \cdots \overset{f_n} \to d_n)$ is an $n$-simplex of $N(D)$ and $x \in F(d_0)_n$. An edge $(f, x)$ is marked in $\gCLM{D}(F)$ if and only if the edge $x \in F(d_0)_1$ is a marked edge.

\begin{rem}
	\label{rel-bar-const}
	The marked simplicial set (total space of) $\gCLM{D}(F)$ is the bar construction $|B(\ast, C, F)|$ of the functor $F:D \to \sSetsM$.
	We recall that $\sSetsMQ$ is a simplicial model category. See \cite{Shul}, \cite{Meyer} for the above notation and a description of a bar construction in a simplicial model category.
\end{rem}

\begin{nota}
	We will abuse notation and denote the total space of $\gCLM{D}(F)$, which is a map of marked simplicial sets having codomain $N(D)$, also by $\gCLM{D}(F)$.
\end{nota}
The following lemma is an easy consequence of the above description of the left Quillen functor $\gCLM{D}$:
\begin{lem}
	\label{fib-img-L}
	For each functor $F:D \to \sSetsMQ$ which is a fibrant object in the projective model category $[D, \sSetsMQ]$, the object $\gCLM{D}(F)$ is of the type $\Nt{\gCLM{D}(F)}$, namely its underlying simplicial map is a coCartesian fibration and the marked edges of the total space are the coCartesian edges. It is therefore a fibrant object of $\sSetsMGen{D}$.
\end{lem}

Based on the above description of the functor $\gCLM{D}$, it is easy to see an inclusion natural transformation:
\begin{equation*}
	\iota_D:\gCLM{D} \Rightarrow \mRNGen{-}{D}.
\end{equation*}
For a functor $F:D \to \sSetsM$ the simplicial map $\iota_D(F)$ determined by the above natural transformation can be described, in degree $n$, as follows:
\[
(\sigma, x) \mapsto (\sigma, x_F),
\]
where $\sigma = (f_1, f_2, \dots, f_n)$ is an $n$-simplex in $N(D) $ and $x_F$ is the natural transformation represented by the following commutative diagram of simplicial sets:
\begin{equation*}
	\xymatrix{
		\Delta[0] \ar[r]^{d^1} \ar[d] & \Delta[1] \ar[r]^{d^2} \ar[d] & \cdots \ar[r]^{d^{n-1}} & \Delta[n-1] \ar[r]^{d^n} \ar[d]_{} & \Delta[n] \ar[d]^{F(f_n \circ \cdots \circ f_1)(x)} \\
		F(\sigma(0)) \ar[r]_{F(f_1)} &  F(\sigma(1)) \ar[r]_{F(f_2)} & \cdots \ar[r]_{F(f_{n-1}) \ \ \ \ \ \ } &  F(\sigma(n-1)) \ar[r]_{ \ \ \ F(f_n)} &  F(\sigma(n))
	}
\end{equation*}
\begin{prop}
	\label{comp-L-Gr-const}
	For each projectively fibrant diagram $F$ in $[D, \sSetsMQ]$, the map $\iota_D(F)$ is a coCartesian equivalence.
\end{prop}
\begin{proof}
	We observe that the map $\iota_D(F)$ induces an isomorphism on the fibers. Now the proposition follows from \cite[Prop. 3.1.3.5.]{JL}
	
\end{proof}

 \begin{coro}
 	\label{homt-func-L}
 	The left Quillen functor $\gCLM{D}$ is a homotopical functor.
 \end{coro}
 \begin{proof}
 	Since $\mRNGen{-}{D}$ is a right Quillen functor, proposition \ref{comp-L-Gr-const} implies that $\gCLM{D}$ preserves weak-equivalences between fibrant objects. Let $H:X \to Y$ be a weak equivalence in the projective model category $[D, \sSetsMQ]$. The chosen fibrant replacement functor $(r, R)$ gives us the following commutative diagram in $[D, \sSetsMQ]$ wherein the horizontal maps are acyclic cofibrations:
 	\begin{equation*}
 		\xymatrix{
 			X \ar[d]_H \ar[r]^{r(X)} & R(X) \ar[d]^{R(H)} \\
 			Y \ar[r]_{r(Y)} & R(Y)
 		}
 	\end{equation*}
 	Applying the left Quillen functor $\gCLM{D}$ to the above commutative square gives us the following commutative square in $\sSetsMGen{D}$:
 	\begin{equation*}
 		\xymatrix@C=16mm{
 			\gCLM{D}(X) \ar[d]_{\gCLM{D}(H)} \ar[r]^{\gCLM{D}(r(X))} & \gCLM{D}(R(X)) \ar[d]^{\gCLM{D}(R(H))} \\
 			\gCLM{D}(Y) \ar[r]_{\gCLM{D}(r(Y))} & \gCLM{D}(R(Y))
 		}
 	\end{equation*}
 	The two horizontal arrows are weak-equivalences because $\gCLM{D}$ preserves acyclic cofibrations. We have observed above that $\gCLM{D}$ preserves weak-equivalences between fibrant objects therefore $\gCLM{D}(R(H)$ is a weak-equivalence. Now the two out of three property of weak-equivalences in a model category implies that $\gCLM{D}(H)$ is a weak-equivalence.
 	
 \end{proof}
 For the terminal map $u:N(D) \to \ast$ induces a \emph{pullback} functor
 \[
 u^*:\sSetsM \to \sSetsMGen{D}
 \]
 In this situation $u^*(X) = (N(D) \times X, p_1)$ for each object $X \in \sSetsM$, where $p_1$ is the projection to $N(D)$. The functor $u^*$ has a left adjoint
 \[
 u_!:\sSetsMGen{D} \to \sSetsM
 \]
 which maps an object $(X, p) \in \sSetsMGen{D}$ to the compopsite $X \to N(D) \to \ast$. It is easy to check that the adjunction $(u_!, u^*)$ is a Quillen adjunction. Now we define a homotopy colimit functor, see definition \ref{hocolim-func}:
 \begin{prop}
 	\label{hocolim-func-mar-SS}
 	The composite right Quillen functor $u_! \circ \gCLM{D}$ is a homotopy colimit functor.
 \end{prop}
 \begin{proof}
 	We have seen above (cor. \ref{homt-func-L}) that the functor $\gCLM{D}$ is homotopical. The functor $u_!$ is also homotopical because it is a left Quillen functor wherein every object of the domain model category is cofibrant.
 	
 	Let $(R, r)$ be the fibrant replacement functor determined by the chosen functorial factorization on the model category $ \sSetsMQ$.
 	We now construct a natural weak equivalence
 	\begin{equation*}
 		\label{hocol-nat-WE}
 		\delta:\mRNGenR{D}u^*R \Rightarrow \Delta R
 	\end{equation*}
 	For a marked simplicial set $X$, the functor $\mRNGenR{D}u^*R(X)$ is defined as follows:
 	\[
 	\mRNGenR{D}u^*R(X)(d) = \mRNGenR{D}(\Sh{N(D)} \times R(X))(d) = \mapFl{\Sh{N(\ovCatGen{d}{D})}}{\Sh{N(D)} \times R(X)}_D
 	\]
 	The functor $\Delta R$ is isomorphic to the functor 
 	\[
 	\mapFl{(\Sh{\Delta[0]}, d)}{\Sh{N(D)} \times R(-)}_D: \sSetsM \to [D, \sSetsM],
 	\]
 	where $(\Sh{\Delta[0]}, d)$ denotes the map $d:\Sh{\Delta[0]} \to \Sh{N(D)}$. Now the desired natural weak-equivalence $\delta$ is defined as follows for each pair $(X, d) \in Ob(\sSetsM) \times Ob(D)$:
 	\[
 	\delta_X(d) := \mapFl{\iota_d^+}{\Sh{N(D)} \times R(X)}_D:\mRNGenR{D}u^*R(X)(d) \to R(X),
 	\]
 	where the map $\iota_d^+$ is defined in \eqref{main-gen-local-mar}. The natural transformation $\delta_X(d)$ is a natural weak equivalence because $\iota_d^+$ is a weak equivalence.
 	
 	\begin{sloppypar}
 		Now it follows from \cite[33.9 (ii)]{DHKS} that the natural weak-equivalence $\mapFl{\iota_d^+}{N(D) \times R(-)}_D$ induces a natural isomorphism between the (derived) functors $Ho(\mRNGenR{D}u^*R)$ and $Ho(\Delta R)$. Since the functor $\Delta$ is homotopical, it follows from \cite[33.9 (ii)]{DHKS} that there is a natural isomorphism between $Ho(\Delta R)$ and $Ho(\Delta)$.
 	\end{sloppypar}
 \end{proof}

The marked version of our homotopy colimit theorem for quasi-categories, namely theorem \ref{Ho-Colim-Thm-QCat-M}, now follows from the above result and proposition \ref{comp-L-Gr-const}.

 Now we state and prove the second prominent theorem of this paper. We closely follow the approach of \cite{HM} in our proof, however in our proof we establish a dual statement to what is proved in the aforementioned paper and thereby we present a new argument which can be suitable adapted to give another proof of \cite[Thm. C]{HM}.
 The theorem will imply that the right derived functor of $\mRNGenR{D}$ and $\mRNGen{-}{D}$ are mutual inverses.
  The above lemma \ref{unit-Eq-mar} will be instrumental in proving this result:
 \begin{thm}
 	\label{main-res}
 	The two  Quillen pairs (left adjoints on top)
 	\begin{equation*}
 		[D, \sSetsMQ] \underset{\mRNGenR{D}} {\overset{\gCLM{}}\rightleftarrows} \sSetsMGen{D} \underset{\mRNGen{-}{D}}{\overset{\pFGen{\bullet}{D}} \rightleftarrows} [D, \sSetsMQ]
 		\end{equation*}
 	 are Quillen equivalences between the coCartesian model category structure on $\sSetsMGen{D}$ and projective model category structure on $[D, \sSetsMQ]$.
 \end{thm}
\begin{proof}
	We will prove this theorem by showing that for each projectively fibrant diagram $F:D \to \sSetsMQ$, the natural transformation $\eta^+_F:F \Rightarrow \mRNGenR{D}\left( \mRNGen{F}{D} \right)$ is objectwise a weak-equivalence in $\sSetsMQ$ and by exhibiting a zig-zag of weak-equivalences in $\sSetsMGen{D}$:
	\begin{equation*}
		\mRNGen{\mRNGenR{D}(X)}{D} \overset{w_1(X)}\leftarrow \gCLM{}\left(\mRNGenR{D}(X)\right) \overset{w_2(X)}\rightarrow X,
		\end{equation*}
	for each fibrant $X \in \sSetsMGen{D}$.
	The first statement follows from lemma \ref{unit-Eq-mar}. Since $X$ is fibrant by assumption, so is $\mRNGenR{D}(X)$ by lemma \ref{Qu-Eq-L-R-mar}. Now proposition \ref{comp-L-Gr-const} tells us that $w_1(X)$ is a weak-equivalence. Finally, we define the map $w_2$ to be the counit map of the adjunction $(\gCLM{} , \mRNGenR{D})$. The component at $X$ of this counit natural transformation, namely $w_2(X)$, is defined (in degree $n$) as follows:
	\begin{equation*}
		(\sigma, x) \mapsto x((\cnGen{n}{\sigma(0)}, id_n)),
		\end{equation*}
	where $\sigma$ is an $n$-simplex in $N(D)$, $x$ is an $n$-simplex in $\mRNGenR{D}(X)(\sigma(0))$ which is a map
	\begin{equation*}
		\xymatrix{
		\Sh{N(\sigma(0)/D)} \times \Fl{\Delta[n]} \ar[rr] \ar[rd] && X \ar[ld] \\
		& \Sh{N(D)}
	  }
		\end{equation*}
	 in $\sSetsMGen{D}$. It is easy to check that both $w_1$ and $w_2$ are natural in $X$.
	The assumption that $X$ is fibrant implies that both the domain and codomain of this counit map $w_2(X)$ are coCartesian fibrations and it is easy to see that $w_2(X)$ induces an equivalence on fibers, therefore it is a weak-equivalence in $\sSetsMGen{D}$. 
	
	The above natural weak-equivalences imply that the right derived functors of $\mRNGenR{D}$ and $\mRNGen{-}{D}$ are mutual inverses of one another. Therefore both right Quillen functors $\mRNGenR{D}$ and $\mRNGen{-}{D}$ are part of Quillen equivalences.
	\end{proof}

The prominent theorem of the paper stated above has the following corollary:
\begin{coro}
	\label{der-funct-L-der-GC-eq}
	The (total) left derived functor of $\gCLM{D}$ is naturally isomorphic to the (total) right derived functor of $\mRNGen{-}{D}$.
	\end{coro}
\begin{nota}
	The total (total) right derived functor of $\mRNGen{-}{D}$ refers to the total right derived functor of the relative nerve functor for marked simplicial sets, see \cite[Prop. 3.2.5.18(2)]{JL}.
	\end{nota}

\section[A homotopy colimit functor]{A homotopy colimit functor for diagrams of quasi-categories}
\label{hoColim-QCat}

The homotopy colimit of a functor $F:D \to \sSetsK$ has a standard construction, namely it is the \emph{diagonal} of the bisimplicial set obtained by applying the nerve functor to the  transport category (Grothendieck construction) of the $n$-simplex functors $F_n$, see \cite[Ch. IV]{GJ}.
In this section we will present a homotopy colimit construction for functors taking values in the Joyal model category of simplicial sets $\sSetsQ$. Our construction is a modification of the aforementioned construction. Our construction will use the homotopy colimit functor constructed on the category of functors taking values in $\sSetsMQ$, see proposition \ref{hocolim-func-mar-SS}. We exhibit a very natural embedding of simplicial sets into marked simplicial sets by marking \emph{equivalences} in simplicial sets. We go on to show that this embedding is an equivalence of the underlying homotopy theories of the two model categories in context. Now the homotopy colimit of a functor $F:D \to \sSetsQ$ is obtained by first composing the functor with the aforementioned embedding $F:D \to \sSetsQ \to \sSetsMQ$, then applying the homotopy colimit functor $u_! \circ \gCLM{D}$ to $F$ and finally inverting all the coCartesian edges of $\gCLM{D}(F) $.

 We begin by recalling the notion of a homotopy colimit functor, see \cite{DS95}, \cite{DHKS} for more detail:
\begin{df}
	\label{hocolim-func}
	A \emph{homotopy colimit functor} on the functor category $[D, \M]$, where $\M$ is a model category, is a homotopical functor
	\[
	hocolim:[D, \M] \to \M
	\]
	such that its induced functor on the homotopy category $Ho(hocolim)$ is a left adjoint to the functor
	\[
	Ho(\Delta):Ho(\M) \to Ho[D, \M].
	\] 
	\end{df}
 In this paper we restrict to those model categories which induce a projective model category structure on the functor category $[D, \M]$.

We recall that an edge in a quasi-category $X$ is called an equivalence if it determines an isomorphism in the homotopy category $Ho(X)$.
We want to extend this definition to all simplicial-sets. We recall that the nerve functor $N:\sSets \to \Cat$ has a left adjoint which is denoted by $\tau_1$, the unit map of this adjunction is denoted by $\eta:id \Rightarrow N\tau_1$
\begin{df}
	\label{eq-sSet}
An edge $y$ of a simplicial set $S$ is called an \emph{equivalence} in the simplicial-set $S$ if its image $\eta_S(y)$ in the quasi-category $N\tau_1(S)$ is an equivalence.
\end{df}
\begin{rem}
Unwinding the definition of the functor $\tau_1$ we observe that an edge $y:a \to b$ in a simplicial set $S$ is an equivalence if and only if there is another edge $\inv{y}:b \to a$ in $S$ and a pair of $2$-simpleces $\sigma, \beta \in S_2$ such that $d_0(\sigma) = \inv{y}$, $d_0(\beta) = y$, $d_2(\sigma) = y$, $d_2(\beta) = \inv{y}$ and $d_1(\sigma) = \unit{a},  d_1(\beta) = \unit{b}$.
\end{rem}
\begin{prop}
	\label{Smap-Eq-to-Eq}
A morphism $F:S \to T$ of simplicial sets maps an equivalence in $S$ to an equivalence in $T$.
\end{prop}
\begin{proof}
Let $F:S \to T$ be a morphism of simplicial-sets and $y$ be an equivalence in $S$.
The unit map of the adjunction $(\tau_1, N)$ provides us with the following commutative square:
\[
\xymatrix{
S \ar[r]^F \ar[d]_{\eta_S} & T \ar[d]^{\eta_T} \\
N\tau_1(S) \ar[r]_{N\tau_1(F)} & N\tau_1(T)
}
\]
By assumption $\eta_S(y)$ is an equivalence in the quasi-category $N\tau_1(S)$. By the above commutative diagram, it is sufficient to show that $N\tau_1(F)(\eta_S(y))$ is an equivalence in the quasi-category $N\tau_1(T)$. The assumption that $\eta_S(y)$ is an equivalence in the quasi-category $N\tau_1(S)$ implies that $y$ is a representative of an isomorphism in $\tau_1(S)$ and the functor $\tau_1(F)$ maps this isomorphism to an isomorphism in $\tau_1(T)$ which determines an equivalence $N\tau_1(F)(\eta_S(y))$ in $N\tau_1(T)$.
\end{proof}

	There is an inclusion map $i:\Delta[1] \hookrightarrow J$, where $J$ is the nerve of the groupoid generated by a single isomorphism. The following proposition is an easy consequence of the definition of an equivalence:
	\begin{prop}
		\label{ext-inc-J}
	 If an edge $y$ in a simplicial set $S$ is an equivalence in $S$ then the morphism $y:\Delta[1] \to S$ can be extended along the inclusion $i$ i.e there is a (dashed) lifting arrow in the following (solid arrow) diagram:
	\[
	\xymatrix{
	\Delta[1] \ar[r]^y \ar[d]_i & S \\
	J \ar@{-->}[ru]_{u(y)}
}
	\]
	\end{prop}

\begin{df}
Let $(S, \E)$ be a marked simplicial set, then the marked arrows determine a simplicial map $\E \times \Delta[1] \to S$. The \emph{localization} of $(S, \E)$ is a simplicial set $S[\inv{\E}]$ defined by the following pushout diagram:
\begin{equation}
\label{local-sSets}
\xymatrix{
\E \times \Delta[1] \ar[r] \ar[d]_{\E \times i} & S \ar[d] \\
\E \times J \ar[r] & S[\inv{\E}]
}
\end{equation}

\end{df}
A characteristic property of the localization $S[\inv{\E}]$ is that for any simplicial map $G:S \to T$ which maps every marked edge in $\E$ to an equivalence in $T$, there exists a unique simplicial map $U:S[\inv{\E}] \to T$ such that the following diagram commutes:
\[
\xymatrix{
	\E \times \Delta[1] \ar[r] \ar[d]_{\E \times i} & S \ar[d]^p \ar@/^1pc/[rdd]^G \\
	\E \times J \ar[r] \ar@/_1pc/[rrd]_{u(\E)} & S[\inv{\E}] \ar@{-->}[rd]^U \\
	& & T
}
\]
The simplicial map $u(\E)$ in the above diagram is determined by the extension map from proposition \ref{ext-inc-J}.
\begin{nota}
	For a marked simplicial set $(S, \E)$, we denote by $(S[\inv{\E}], p(\E))$ the marked simplicial set whose set of marked arrows is the image of $\E$ under the projection map $p:S \to S[\inv{\E}]$. We observe that the set $p(\E)$ is a subset of the set of equivalences in $S[\inv{\E}]$ and therefore we have an inclusion map $(S[\inv{\E}], p(\E)) \subseteq (S[\inv{\E}], Eq)$.
	\end{nota}
 
 \begin{lem}
 	\label{inc-acy-cof}
 	The inclusion map $\iota_{(S, \E)}:(S[\inv{\E}], p(\E)) \subseteq (S[\inv{\E}], Eq)$ is an acyclic cofibration in $\sSetsQ$.
 	\end{lem}
 \begin{proof}
 	In light of \cite[Prop. 4.22]{sharma} it is sufficient to show that $\iota_{(S, \E)}$ has the left lifting property with respect to all fibrations between fibrant objects in $\sSetsQ$. Let $q:(X, Eq) \to (Y, Eq)$ be such a fibration. We want to show that whenever we have the following (solid) commutative diagram in $\sSetsM$, there exists a lifting arrow:
 		\[
 	\xymatrix{
 		(S[\inv{\E}], p(\E)) \ar[r] \ar[d]_{\iota_{(S, \E)}} & (X, Eq) \ar[d]^{q} \\
 		(S[\inv{\E}], Eq) \ar[r]  & (Y, Eq)
 	}
 	\]
 	We recall that the underlying simplicial map of $\iota_{(S, \E)}$ is the identity map.
 	Applying the forgetful functor $U$ to the above diagram, we get the following (outer) commutative diagram in $\sSets$ wherein the existence of the lifting arrow is obvious:
 	\[
 	\xymatrix{
 	(S[\inv{\E}]) \ar[r] \ar@{=}[d] &X \ar[d]^{U(q)} \\
 	(S[\inv{\E}]) \ar[r] \ar@{-->}[ru] & Y
    }
 	\]
 	By proposition \ref{Smap-Eq-to-Eq}, each simplicial morphism maps equivalences to equivalences which implies that the following (solid) commutative diagram has a lifting arrow:
 	\[
 \xymatrix{
 	(S[\inv{\E}], p(\E)) \ar[r] \ar[d]_{\iota_{(S, \E)}} & (X, Eq) \ar[d]^{q} \\
 	(S[\inv{\E}], Eq) \ar[r] \ar@{-->}[ru] & (Y, Eq)
 }
 \]
 Thus we have shown that the map $\iota_{(S, \E)}$ is an acyclic cofibration in $\sSetsMQ$.
 	\end{proof}

The above localization defines a functor $L:\sSetsM \to \sSets$. This functor has a right adjoint $E:\sSets \to \sSetsM$ which maps a simplicial $S$ to a marked simplicial set $(S, Eq)$ where $Eq$ is the set of equivalences in $S$. For each marked simplicial set $(S, \E)$, the unit map of this adjunction is the following composite:
\[
  (S, \E) \overset{p}\to (S[\inv{\E}], p(\E)) \subseteq (S[\inv{\E}], Eq) = EL((S, \E))
\]  
 where $p$ is the projection map. The counit map is an isomorphism.

\begin{prop}
	For each marked simplicial set $(S, \E)$ the projection map $p:S \to S[\inv{\E}]$ is a cofibration.
	\end{prop}
\begin{proof}
	The inclusion map $i:\Delta[1] \to J$ is a cofibration of simplicial sets. Now the proposition follows from the observation that cofibrations are preserved under cobase changes in a model category.
	\end{proof}
	We observe that the pushout square \eqref{local-sSets}, in $\sSets$, determines the following commutative square in the category of marked simplicial sets:
	\begin{equation}
	\label{mar-local}
	\xymatrix{
		\Sh{\E} \times \Sh{\Delta[1]} \ar[r] \ar[d]_{\E \times \Sh{i}} & (S, \E) \ar[d]^p \\
		\Sh{\E}  \times \Sh{J} \ar[r] & (S[\inv{\E}], p(\E))
	}
	\end{equation}
	The above proposition holds more strongly in the category of marked simplicial sets. We want to show that the unit map of the above adjunction is an acyclic cofibration in $\sSetsMQ$.

\begin{prop}
	For each marked simplicial set $(S, \E)$, the  commutative square \eqref{mar-local}, in $\sSetsM$, is a pushout square
	\end{prop}
\begin{proof}
	Let $(T, \Delta)$ be another marked simplicial set such that we have the following (outer) commutative diagram:
	\begin{equation}
	\label{mar-pushout}
	\xymatrix{
		\Sh{\E} \times \Sh{\Delta[1]} \ar[r] \ar[d]_{\E \times i} & (S, \E) \ar[d]^p \ar@/^1pc/[rdd]^G \\
		\Sh{\E}  \times \Sh{J} \ar[r] \ar@/_1pc/[rrd]_{u(\E)} & (S[\inv{\E}], p(\E)) \ar@{-->}[rd]^L \\
		& & (T, \Delta)
	}
	\end{equation}
	We want to show the existence of the (dotted) lifting arrow $L$ which makes the entire diagram commutative.
	In this situation, the maps $G$ and $u(\E)$ take values in the (marked) simplicial subset $(T, Eq \cap \Delta) \subseteq (T, \Delta)$. Applying the forgetful functor $U$ to \eqref{mar-local} we get the  commutative diagram \eqref{local-sSets} in $\sSets$ which is a pushout diagram, therefore there exists a (dotted) arrow $L$ in the following (outer) commutative diagram in $\sSets$:
	\[
	\xymatrix{
		\E \times \Delta[1] \ar[r] \ar[d]_{\E \times i} & S \ar[d]^{U(p)} \ar@/^1pc/[rdd]^{U(G)} \\
		\E \times J \ar[r] \ar@/_1pc/[rrd]_{U(u(\E))} & S[\inv{\E}] \ar@{-->}[rd]^{L} \\
		& & T
	}
	\]
 The commutativity of this (simplicial) diagram implies that $L$ maps the set of edges $p(\E)$ into $(T, Eq \cap \Delta)$. This implies that $L$ is the desired dotted arrow in the commutative diagram \eqref{mar-pushout} in $\sSetsM$.
	
	\end{proof}
\begin{coro}
	\label{unit-Eq}
	The unit map of the adjunction $(L, E)$ is a (natural) acyclic cofibration.
	\end{coro}
\begin{proof}
	Since $\Sh{\E}$ is a discrete (marked) simplicial set therefore the commutative square \eqref{mar-local} is can be rewritten as follows:
	\begin{equation*}
	\xymatrix{
		\underset{\E} \sqcup \ \Sh{\Delta[1]} \ar[r] \ar[d]_{\underset{\E} \sqcup \ \Sh{i}} & (S, \E) \ar[d]^p \\
		\underset{\E} \sqcup \ \Sh{J} \ar[r] & (S[\inv{\E}], p(\E))
	}
	\end{equation*}
	The morphism $\Sh{i}$ is an acyclic cofibration in $\sSetsMQ$, it follows from\cite[Prop. 7.2.12]{Hirchhorn} that $\underset{\E} \sqcup \ \Sh{i}$ is also an acyclic cofibration in $\sSetsMQ$. The cobase change of an acyclic cofibration is again an acyclic cofibration which means that projection map $p$ in \eqref{mar-local} is an acyclic cofibration.
	Now the corollary follows from lemma \ref{inc-acy-cof}.
	\end{proof}

\begin{prop}
	The adjunction $(L, E)$ induces an (adjoint) equivalence between the homotopy category of  $\sSetsMQ$ and that of $\sSetsQ$.
	\end{prop}
\begin{proof}
	We have seen above that both the unit and the counit maps of the adjunction $(L, E)$ are (natural) weak equivalences therefore, in light of \cite[Prop. ]{DHKS} it is sufficient to show that $L$ and $E$ are homotopical (weak-equivalence preserving) functors. Let $F:S \to T$ be a weak equivalence in $\sSetsQ$. We want to show that $E(F):E(S) \to E(T)$ is a weak equivalence in $\sSetsMQ$. Since the counit map of the adjunction in context is an isomorphism therefore the right adjoint $E$ is fully faithful. We observe that for all $n \in \Nat$, $E(S \times \Delta[n]) = E(S) \times \Fl{\Delta[n]}$. These two facts together imply that for any simplicial set $W$, we have a natural isomorphism:
	\[
	[S, W] \cong \Fl{[(S, Eq),(W, Eq)]} =\Fl{[E(S), E(W)]}.
	\]
	We observe that for any fibrant marked simplicial set $Z$ we have the following equality: $Z = EU(Z)$, where $U$ is the forgetful functor. Now the following commutative diagram implies that $E$ preserves weak-equivalences:
	\begin{equation*}
	\xymatrix{
	 [S, U(Z)] \ar[r]  & \Fl{[E(S), EU(Z)]}  \ar@{=}[r] & \Fl{[E(S), Z]}  \\
	 [T, U(Z)] \ar[r] \ar[u]^{[F, U(Z)]}  & \Fl{[E(T), EU(Z)]} \ar@{=}[r] \ar[u] & \Fl{[E(T), Z]} \ar[u]_{[E{F}, Z]}
    }
	\end{equation*}
	
	Let $G:(S, \E) \to (T, \Delta)$ be a weak equivalence in $\sSetsMQ$. The unit natural transformation $\eta$ provides us the following commutative diagram:
	\[
	\xymatrix{
	(S, \E) \ar[r]^G \ar[d]_{\eta_{(S, \E)}} & (T, \Delta) \ar[d]^{\eta_{(T, \Delta)}} \\
	EL((S, \E)) \ar[r]_{EL(G)} & EL((T, \Delta))
    }
	\]
	The vertical arrows are acyclic cofibrations by corollary \ref{unit-Eq} and the top horizontal arrow is a weak equivalence by assumption, therefore the two out of three property of weak-equivalences in a model category implies that $EL(G)$ is a weak-equivalence in $\sSetsMQ$. Let $Y$ be any quasi-category, then we have the following commutative square:
	\begin{equation*}
	\xymatrix{
		[L((S, \E)), Z] \ar[r]^\cong   & \Fl{[EL((S, \E)), E(Z)]}  \\
		[L((T, \Delta)), Z] \ar[r]_\cong \ar[u]^{[L(G), Z]}  & \Fl{[EL((T, \Delta)), E(Z)]} \ar[u]_{[EL(G), E(Z)]}
	}
	\end{equation*}
	The two horizontal arrows are isomorphisms and the right downward arrow is a weak-equivalence. The two out of three property of weak equivalences in a model category implies that for each quasi-category $Z$, the simplicial map $[L(G), Z]$ is a weak equivalence and therefore $L(G)$ is a weak equivalence in $\sSetsQ$. Thus we have shown that $L$ preserves weak equivalences.
	
	\end{proof}
	 For any (small) category $D$, the adjunction $(L, E)$ induces an adjunction 
	 \begin{equation*}
	 L^D:[D, \sSets^+] \rightleftharpoons [D, \sSets]:E^D.
	 \end{equation*}
	 Each of the two functor categories can be endowed with a \emph{projective} model category structure which is inherited from $\sSetsMQ$ and $\sSetsQ$ respectively, see \cite[Remark 2.8.6]{JL}.
	 \begin{nota}
	 	We will denote the aforementioned model categories by $[D, \sSetsMQ]$ and $[D, \sSetsQ]$ respectively.
	 	\end{nota}
 	\begin{nota}
 		We will denote the homotopy categories of $[D, \sSetsMQ]$ and $[D, \sSetsQ]$ by $Ho(\sSetsMQ^D)$ and $Ho(\sSetsQ^D)$ respectively.
 		\end{nota}
	\begin{coro}
	The adjunction $(L^D, E^D)$ induces an (adjoint) equivalence of categories on the homotopy categories:
	\begin{equation*}
	Ho(L^D):Ho(\sSetsMQ^D) \rightleftharpoons Ho(\sSetsQ^D):Ho(E^D).
	\end{equation*}
	\end{coro}
\begin{proof}
	Both adjoint functors $L^D$ and $E^D$ are homotopical functors and therefore both are deformable in the sense of \cite[Def. ]{DHKS} with identity deformations. Now the result follows from \cite[Prop. 45.2]{DHKS}.

	\end{proof}
 We denote the colimit functor on the functor category $[D, \sSetsM]$ by
 \[
 {\varinjlim}^+:[D, \sSetsM] \to \sSetsM
 \]
 Consider the following composite functor:
 \begin{equation}
 \label{Eq-colim}
 [D, \sSets] \overset{E^D} \to [D, \sSetsM] \overset{\colim{+}} \to \sSetsM \overset{L} \to \sSets
 \end{equation}
 which we denote by $L\colim{}$. The next lemma shows that this functor is a colimit functor on $\sSets$:
 \begin{lem}
 	\label{colim-in-QCat}
 	The composite functor $L\colim{}$, \eqref{Eq-colim}
 	is a colimit functor.
 	\end{lem}
 \begin{proof}
 	Let $F:D \to \sSetsQ$ be a functor and $Y$ be a simplicial set.
 	It is sufficient to observe the following chain of (natural) bijections:
 	\begin{multline*}
 	\sSets(L(\colim{}(E^D(F))), Y) \cong \sSetsM(\colim{+}(E^D(F)), E(Y)) \cong [D, \sSetsM](E^D(F), \Delta(E(Y))) = \\
 	[D, \sSetsM](E^D(F), E^D(\Delta(Y))) \cong [D, \sSets](L^D(E^D(F)), \Delta(Y)) \cong \\ \sSets(\colim{}(L^D(E^D(F)), Y) \cong \sSets(\colim{}(F), Y).
 	\end{multline*}
 	The last bijection in the above chain follows from the observation that the counit of the adjunction $(L^D, E^D)$ is an isomorphism.
 
 	\end{proof}
 We now consider the following composite functor:
 \begin{equation}
 \label{Eq-hocolim}
 [D, \sSets] \overset{E^D} \to [D, \sSetsM] \overset{u_! \circ \gCLM{D}} \to \sSetsM \overset{L} \to \sSets
 \end{equation}
 which we denote by $L\hColim{}$.
 Now we are ready to state and prove the main result of this section:
 \begin{thm}
 	The functor $L\hColim{}$ is a homotopy colimit functor.
 	\end{thm}
 \begin{proof}
 	The functor $L\hColim{}$ is a composite of three homotopical functors therefore it is homotopical. Now we consider the following chain of (natural) bijections for a pair of objects $F \in [D, \sSets]$ and $Y \in \sSets$:
 	\begin{multline*}
 	Ho_{\sSets}(L\hColim{}(F), Y) \cong Ho_{\sSetsM}(u_!\gCLM{D}E^D(F), E(Y)) \cong \\
 	Ho_{\sSetsMQ^D}(E^D(F), \Delta(E(Y))) = Ho_{\sSetsMQ^D}(E^D(F), E^D(\Delta(Y))) \cong \\
 	Ho_{\sSetsQ^D}(L^DE^D(F), \Delta(Y)) \cong Ho_{\sSetsQ^D}(F, \Delta(Y))
 	\end{multline*}

 	\end{proof}
 
 Our homotopy colimit theorem for quasi-categories, namely theorem \ref{Ho-Colim-Thm-QCat}, now follows from the above result and proposition \ref{comp-L-Gr-const}.

 \appendix
 \section{The classical Grothendieck construction}
\label{Classic-GC}

The classical Grothendieck construction \cite{SGA} associates to a pseudo-functor $F:C \to \Cat$, an opfibration over the category $C$ which we denote as follows:
\[
p:\rNGen{F}{C} \to C.
\]
Throughout this paper we will abuse notation and refer to the domain category $\rNGen{F}{C}$ of the opfibration $p$ (also known as the total category of $p$) as the Grothendieck construction of $F$.
This assignment defines a bifunctor
\begin{equation*}
\rNGen{-}{C}:\CatHom{C}{\Cat}{Ps} \to \Cat/C.
\end{equation*}
This bifunctor factors through a bicategory of opfibrations $OpFib(C)$ over $C$ as follows:
\begin{equation*}
\CatHom{C}{\Cat}{Ps} \overset{\rNGen{-}{C}} \to OpFib(C) \hookrightarrow \Cat/C.
\end{equation*}
The first arrow in the above diagram is an equivalence of bicategories. We refer the interested reader to the book \cite{PJ} for further details.

 A diagram in $\Cat$ is a functor from a small category into $\Cat$.
We begin this section by a review of the classical Grothendieck construction for a diagram in $\Cat$. The main objective of this section is to associate a \emph{double category} (\emph{i.e.} a category object in $\Cat$) to a diagram $F:C \to \Cat$ and obtain the Grothendieck construction $\rNGen{F}{C}$ by passing to the object sets and object functions of all categories and functors involved in the definition of the aforementioned category object associated with $F$. The category obtained this way from a double category is often called the \emph{horizontal structure} of the double category.

 The object set of $\rNGen{F}{C}$ is the following:
\[
Ob\left( \rNGen{F}{C} \right) = \underset{c \in C} \bigsqcup Ob(F(c)).
\]
In order to describe the morphisms we consider the (categorical) \emph{mapping path space} construction of a functor $F(f):F(c_1) \to F(c_2)$, where  $f:c_1 \to c_2$ is a morphism in $C$, which is described by the following cartesian diagram in $\Cat$:
\begin{equation*}
\xymatrix{
\P^1_{F(f)} \ar[r] \ar[d] & F(c_2)^I \ar[d]^{F(c_2)^{i_0}} \\
F(c_1) \ar[r]_{F(f)} & F(c_2)
}
\end{equation*}
where $i_0:0 \to I$ is the inclusion into the \emph{source}. The category $\P^1_{F(f)}$ can be described as follows: An object in $\P^1_{F(f)}$ is a pair
\[
 (\beta_0, \beta_1) \in F(s(f)) \times F(t(f))^I.
 \]
 A morphism in $\P^1_{F(f)}$, from $(\beta_0, \beta_1)$ to $(\gamma_0, \gamma_1)$ is a pair of maps 
 \[
 (f_0, f_1) \in F(s(f)) \times F(t(f)),
 \]
 such that the following diagram commutes:
 \begin{equation*}
 \xymatrix{
F(f)(\beta_0) \ar[d]_{F(f)(f_0)} \ar[r]^{\beta_1} & \beta_1(1) \ar[d]^{f_1} \\
F(f)(\gamma_0) \ar[r]_{\gamma_1} & \gamma_1(1)
}
 \end{equation*}

 Now we can describe the morphism set of $\rNGen{F}{C}$ as follows:
\[
Mor\left(\rNGen{F}{C}\right) = \underset{f \in Mor(C)} \bigsqcup Ob(\P^1_{F(f)}).
\]

We observe that category $\P^1_{F(f)}$ is equipped with the following two functors which we refer to as the source and target functors:
\begin{equation*}
\label{src-MPC}
s:\P^1_{F(f)} \to F(s(f)) \hookrightarrow \underset{c \in C} \bigsqcup F(c).
\end{equation*}
This functor is defined on objects as follows:
\[
s((\beta_0, \beta_1)) := \beta_0
\]
and it is defined on morphisms as follows:
\[
s((f_0, f_1)) := f_0.
\]
The second functor is 
\begin{equation*}
\label{tar-MPC}
t:\P^1_{F(f)} \to F(t(f)) \hookrightarrow \underset{c \in C} \bigsqcup F(c).
\end{equation*}
This functor is defined on objects as follows:
\[
t((\beta_0, \beta_1)) :=  \beta_1(1)
\]
and it is defined on morphisms as follows:
\[
t((f_0, f_1)) :=  f_1.
\]
The source and target functors described above, define the following graph in the category $\Cat$:
\begin{equation*}
\underset{f \in C} \bigsqcup \P^1_{F(f)} \overset{s} {\underset{t} \rightrightarrows}  \underset{c \in C} \bigsqcup F(c).
\end{equation*}
We claim that the above graph is the underlying graph of a category object in $\Cat$. In order to construct this category object, we will have to define two more functors, the first of these functors is a \emph{unit} functor 
\[
u:\underset{c \in C} \bigsqcup F(c) \to \underset{f \in C} \bigsqcup \P^1_{F(f)}.
\]
In order to do so, it is sufficient to define a functor
\[
u_c:F(c) \to \P^1_{F(id_c)},
\]
for each $c \in C$. This functor
is defined on objects as follows:
\[
u_c(\beta_0) := (\beta_0, id_{F(\beta_0)}).
\]
It is defined on morphisms as follows:
\[
u_c(f_0) := (f_0, f_0).
\]
The second functor needed to define the category object in $\Cat$ is the \emph{composition} functor:
\begin{equation*}
-\circ-:\underset{f \in C} \bigsqcup \P^1_{F(f)} \underset{s = t} \times \underset{f \in C} \bigsqcup \P^1_{F(f)} \to \underset{f \in C} \bigsqcup \P^1_{F(f)}.
\end{equation*}
In light of the following isomorphism of categories:
\begin{equation*}
\underset{f \in C} \bigsqcup \P^1_{F(f)} \underset{s = t} \times \underset{f \in C} \bigsqcup \P^1_{F(f)} \cong  \underset{(f_1, f_2) \in N(C)_2} \bigsqcup \P^1_{F(f_2)} \underset{s = t}\times \P^1_{F(f_1)}
\end{equation*}

 it is sufficient to define a functor
\[
-\circ-:\P^1_{F(f_2)} \underset{s = t}\times \P^1_{F(f_1)} \to \P^1_{F(f_2f_1)},
\]
for each pair of composable arrows $(f_1, f_2)$ in $C$.

Let $(\beta_0, \beta_1) \in Ob(\P^1_{F(f_1)})$ and $(\gamma_0, \gamma_1) \in Ob(\P^1_{F(f_2)})$ be two objects in $\underset{f \in C} \bigsqcup \P^1_{F(f)}$ such that
\[
s((\gamma_0, \gamma_1)) = t((\beta_0, \beta_1)).
\]
Now we define
\[
(\gamma_0, \gamma_1) \circ (\beta_0, \beta_1) := (\beta_0, \gamma_1 \circ F(f_2)(\beta_1)).
\]
Let $(g_0, g_1) \in Mor(\P^1_{F(f_1)})$ and $(h_0, h_1) \in Ob(\P^1_{F(f_2)})$ be two objects in $\underset{f \in C} \bigsqcup \P^1_{F(f)}$ such that
\[
s((h_0, h_1)) = t((g_0, g_1)).
\]
Now we define
\[
(h_0, h_1) \circ (g_0, g_1) := (g_0, h_1).
\]

One can now check that the functors $(s, t, u, - \circ -)$ associate to the functor $F:C \to \Cat$, the following category object in $\Cat$:
\begin{equation*}
\left(\underset{c \in C} \bigsqcup F(c), \underset{f \in C} \bigsqcup \P^1_{F(f)}, s, t, u, -\circ-  \right).
\end{equation*}
The Grothendieck construction of the functor $F$, namely the category $\rNGen{F}{C}$ is now obtained by passing to the object sets and object functions of all categories and functors involved in the definition of the above category object as follows:
\begin{equation*}
\left(Ob\left(\underset{c \in C} \bigsqcup F(c) \right), Ob\left(\underset{f \in C} \bigsqcup \P^1_{F(f)} \right), Ob(s), Ob(t), Ob(u), Ob(-\circ-)  \right).
\end{equation*}

A category object in $\Cat$ is also known as a \emph{double category}. A category obtained from a double category by passing to the object sets and object functions of all categories and functors defining the double category, like above, is called the \emph{horizontal structure} of the double category.

 \section[Marked simplicial sets]{A review of marked simplicial sets}
\label{mar-sSets}
In this appendix we will review the theory of marked simplicial sets. Later in this paper we will develop a theory of coherently commutative monoidal objects in the category of marked simplicial sets.

\begin{df}
	\label{mar-sSet}
	A \emph{marked} simplicial set is a pair $(X, \E)$, where $X$ is a simplicial set and $\E$ is a set of edges of $X$ which contains every degenerate edge of $X$. We will say that an edge of $X$ is \emph{marked} if it belongs to $\E$.
	A morphism $f:(X, \E) \to (X', \E')$ of marked simplicial sets is a simplicial map $f:X \to X'$ having the property that $f(\E) \subseteq \E'$. We denote the category of marked simplicial sets by $\sSetsM$.
	\end{df}

Every simplicial set $S$ may be regarded as a marked simplicial set in many ways. We mention two extreme cases: We let $\Sh{S} = (S, S_1)$ denote the marked simplicial set in which every edge is marked. We denote by $\Fl{S} = (S, s_0(S_0))$ denote the marked simplicial set in which only the degenerate edges of $S$ have been marked.

The category $\sSetsM$ is \emph{cartesian-closed}, \emph{i.e.} for each pair of objects $X, Y \in Ob(\sSetsM)$, there is an internal mapping object $[X, Y]^+$ equipped with an \emph{evaluation map} $[X, Y]^+ \times X \to Y$ which induces a bijection:
\[
\sSetsM(Z, [X, Y]^+) \overset{\cong} \to \sSetsM(Z \times X, Y),
\]
for every $Z \in \sSetsM$.
\begin{nota}
	We denote by $\Fl{[X, Y]}$ the underlying simplicial set of $[X, Y]^+$.
	\end{nota}
The mapping space $\Fl{[X, Y]}$ is characterized by the following bijection:
\[
\sSets(K, \Fl{[X, Y]}) \overset{\cong} \to \sSetsM(\Fl{K} \times X, Y),
\]
for each simplicial set $K$.
\begin{nota}
	We denote by $\Sh{[X, Y]}$ the simplicial subset of $\Fl{[X, Y]}$ consisting of all simplices $\sigma \in \Fl{[X, Y]}$ such that every edge of $\sigma$ is a marked edge of $[X, Y]^+$.
\end{nota}
The mapping space $\Sh{[X, Y]}$ is characterized by the following bijection:
\[
\sSets(K, \Sh{[X, Y]}) \overset{\cong} \to \sSetsM(\Sh{K} \times X, Y),
\]
for each simplicial set $K$.

 The Joyal model category structure on $\sSets$ has the following analog for marked simplicial sets:
 \begin{thm}
 	\label{Joyal-sSetsM}
 	There is a left-proper, combinatorial model category structure on the category of marked simplicial sets $\sSetsM$ in which a morphism $p:X \to Y$ is a
 	\begin{enumerate}
 		\item cofibration if the simplicial map between the underlying simplicial sets is a cofibration in $\sSetsQ$, namely a monomorphism.
 		
 		\item a weak-equivalence if the induced simplicial map on the mapping spaces
 		\[
 		\Fl{[p, \Nt{K}]}:\Fl{[X, \Nt{K}]} \to \Fl{[Y, \Nt{K}]}
 		\]
 		is a weak-categorical equivalence, for each quasi-category $K$.
 		
 		\item fibration if it has the right lifting property with respect to all maps in $\sSetsM$ which are simultaneously cofibrations and weak equivalences.
 		
 		\end{enumerate}
 	Further, the above model category structure is enriched over the Joyal model category, i.e. it is a $\sSetsQ$-model category.
 	\end{thm}
 The above theorem follows from \cite[Prop. 3.1.3.7]{JL}.
 \begin{nota}
 	We will denote the model category structure in Theorem \ref{Joyal-sSetsM} by $\sSetsMQ$ and refer to it either as the \emph{Joyal} model category of \emph{marked} simplicial sets or as the model category of marked quasi-categories.
 	\end{nota}
 \begin{thm}
 	\label{Cart-cl-Mdl-S-plus}
 	The model category $\sSetsMQ$ is a cartesian closed model category.
 	\end{thm}
 \begin{proof}
 	The theorem follows from \cite[Corollary 3.1.4.3]{JL} by taking $S = T = \Delta[0]$.
 	\end{proof}
 
 There is an obvious forgetful functor $U:\sSetsM \to \sSets$. This forgetful functor has a left adjoint $\Fl{(-)}:\sSets \to \sSetsM$.
 \begin{thm}
 	\label{Quil-eq-JQ-MJQ}
 	The adjoint pair of functors $(\Fl{(-)}, U)$ determine a Quillen equivalence between the Joyal model category of marked simplicial sets and the Joyal model category of simplicial sets.
 	\end{thm}
 The proof of the above theorem follows from \cite[Prop. 3.1.5.3]{JL}.
 \begin{rem}
 	A marked simplicial set $X$ is fibrant in $\sSetsMQ$ if and only if it is a quasi-category with the set of all its equivalences as the set of marked edges.
 	\end{rem}
 The following proposition brings out a very important distinction between the model category $\sSetsMQ$ and $\sSetsQ$:
 \begin{prop}
 	\label{inc-J-acy-cof}
 	The inclusion map
 	\[
 	\Sh{i}:\Sh{\Delta[1]} \to \Sh{J}
 	\]
 	is an acyclic cofibration in the model category $\sSetsMQ$.
 	\end{prop}
 \begin{proof}
 	In light of \cite[Prop. 4.22]{sharma} it is sufficient to show that $\Sh{i}$ has the left lifting property with respect to every fibration between fibrant objects in $\sSetsMQ$. Let $p:X \to Y$ be such a fibration. We observe that $p = EU(p)$. By adjointness the aforementioned left lifting property is equivalent to $U(p)$ having left lifting property with respect to $L(\Sh{i})$. The proposition now follows from the observation that the map $L(\Sh{i})$ is a canonical isomorphism in $\sSetsQ$.
 	\end{proof}
 The above proposition does NOT hold in $\sSetsQ$ \emph{i.e.} the inclusion map $i:\Delta[1] \to J$ is a cofibration but it is NOT acyclic in $\sSetsQ$.
  \section[Comparison with Relative Nerve]{Comparison with Relative Nerve}
\label{comp-rel-ner}
The notion of \emph{relative nerve} was introduced in 
\cite[Sec. 3.2.5]{JL}. In this appendix we will show that our definition of Lurie's Grothendieck construction of a diagram of simplicial sets, namely definition \ref{LGrC}, is essentially the same as that of the relative nerve of the diagram \cite[Defn. 3.2.5.2.]{JL}. More precisely, we show that for any diagram of simplicial sets $F:D ]to \sSets$, we have the following (natural) isomorphism:
\begin{equation*}
	\rNGen{F}{D} \cong N_F(D),
	\end{equation*}
where $N_F(D)$ denotes the relative nerve of $F$.
 We begin by reviewing the relative nerve:
\begin{df}
	\label{rel-Ner}
	Let $D$ be a category, and $F:D \to \sSets$ a functor. The nerve of $D$ relative to $F$ is the simplicial set $N_F (D)$ whose $n$-simplices are sets consisting of:
	\begin{enumerate}
		\item[(i)] a functor $d:[n] \to C$; We write $d(i, j)$ for the image of $i \le j$ in $[n]$.
		\item[(ii)] for every nonempty subposet $J \subseteq [n]$ with maximal element $j$, a map $\tau^J:\Delta^J \to F(d(j))$,
		\item[(iii)] such that for nonempty subsets $I \subseteq J \subseteq [n]$ with respective maximal elements  $i \le j$, the following diagram commutes:
		\[
		\xymatrix{
			\Delta^I\ar[r]^{\tau^I} \ar@{_{(}->}[d] & F(d(i)) \ar[d]^{F(d(i, j))} \\
			\Delta^J\ar[r]_{\tau^J}  & F(d(j))
		}
		\]
	\end{enumerate}
\end{df}

For any functor $F:D \to \sSets$, there is a canonical map $p_F: N_F(D) \to N(D)$ to the ordinary nerve of $D$, induced by the unique map to the terminal object $\Delta^0 \in \sSets$ \cite[ 3.2.5.4]{JL}. When $F$ takes values in quasi-categories, this canonical map is a coCartesian fibration.
\begin{rem}
	\label{Rel-Ner-edge}
	A vertex of the simplicial set $N_F(D)$ is a pair $(c, g)$, where $c \in Ob(D)$ and $g \in F(D)_0$.
	An edge $\ud{e}:(c, g) \to (d, k)$ of the simplicial set $N_F $ consists of a pair $(e, h)$, where $e:D \to d$ is an arrow in $D$ and $h:f(e)_0(g) \to k$ is an edge of $F(d)$.
\end{rem}
An immediate consequence of the above definition is the following proposition:
\begin{prop}
	\label{Rel-Ner-isom-func-val}
	Let $F:D \to \sSets$ be a functor, then the fiber of $p_F:N_F(D) \to N(D)$ over any $c \in Ob(D)$ is isomorphic to the simplicial set $F(D)$.
\end{prop}

The following lemma is a consequence of this definition and the above discussion:

\begin{lem}
	\label{isom-Rel-Ner}
	For each functor $X:D \to \sSets$, we have the following isomorphism in the category $\ovCatGen{\sSets}{N(D)}$:
	\[
	\rNGen{X}{C} \cong N_X(D).
	\]
\end{lem}
\begin{proof}
	An $n$-simplex in $\rNGen{X}{C}$ is a pair $(\sigma, \beta)$, where $\sigma \in N(D)_n$ and $\beta:D(n) \Rightarrow X(\sigma)$ is a natural transformation, see Definition \ref{nth-can-func-Path-Fib}. The $n$-simplex $\sigma$ can be described as a functor $\sigma:[n] \to C$. The inclusion of each non-empty subposet $i_J:J \subseteq [n]$ gives a map
	\[
	\left(\rNGen{X}{C} \right)(i_J):\left(\rNGen{X}{C} \right)_n \to \left(\rNGen{X}{C} \right)_J.
	\]
	We are using the fact that $J$ is isomorphic to an object of $\Delta$ which we also denote by $J$. The inclusion map can now be seen as a map in $\Delta$.
	This map gives us a $J$-simplex
	\[
	\left(\rNGen{X}{C} \right)(i_J)((\sigma, \beta)) := (\sigma^\ast, \beta^\ast),
	\]
	where $\sigma^\ast$ is the composite
	
	\[
	J \hookrightarrow [n] \overset{\sigma} \to C
	\]
	and $\beta^\ast$ is the following composite:
	\begin{equation*}
	\xymatrix{
		&  & \ar@{=>}[dd]^{\beta} &\\
		J \ar@{^{(}->}[r]^{} & [n] \ar@/^3pc/[rr]^{c(n)} \ar@/_3pc/[rr]_{X(\sigma)} && \sSets \\
		& &
	}
	\end{equation*}
	
	Now the above natural transformation $\beta^\ast$ gives us a simplicial map:
	\[
	\beta^\ast(j'):\Delta^J \to X(\sigma(j'))
	\]
	where $j'$ is the maximal element of $J$. For an inclusion $J' \subseteq J$, condition $(iii)$ of definition \ref{rel-Ner} is satisfied because the composite $J' \subseteq J \subseteq [n]$ determines a composite map in $\Delta$. This defines a map $f_n:\left(\rNGen{X}{C} \right)_n \to N_X(D)_n$ for all $n \ge 0$. One can check that this collection glues into a simplicial map.
	
	Now we define the inverse map. An $n$-simplex $\gamma$ in $N_X(D)$ contains a functor $\sigma:[n] \to C$, by \ref{rel-Ner}$(i)$, which uniquely determines an $n$-simplex $\sigma$ of $N(D)$. We recall that an $n$-simplex in $\rNGen{X}{C} $ is a pair $(\sigma, \tau)$, where $\tau:D(n) \Rightarrow X(\sigma)$ is a natural transformation.
	We observe the following sequence of inclusion maps in the category $\Delta$:
	\begin{equation*}
	[0] \subset [1] \subset \cdots \subset [i] \subset \cdots \subset [n].
	\end{equation*}
	Now, it follows from \ref{rel-Ner}$(ii)$ and \ref{rel-Ner}$(iii)$ that the above sequence of inclusions gives us the following commutative diagram:
	\begin{equation*}
	\xymatrix{
		\Delta[0] \ar[r] \ar[d]^{\tau^{[0]}} & \Delta[1] \ar[r] \ar[d]^{\tau^{[1]}} & \cdots \ar[r] & \Delta[n-1] \ar[r] \ar[d]^{\tau^{[n-1]}} & \Delta[n] \ar[d]^{\tau^{[n]}} \\
		X(\sigma(0)) \ar[r] & X(\sigma(1)) \ar[r] & \cdots \ar[r] & X(\sigma(n-1)) \ar[r] & X(\sigma(n))
	}
	\end{equation*}
	which defines a natural transformation $\tau:D(n) \Rightarrow X(\sigma)$. This defines a map
	\[
	\inv{f}_n:N_X(D)_n \to \left( \rNGen{X}{C} \right)_n
	\]
	which maps $\gamma$ to $(\sigma, \tau)$, for all $n \ge 0$.
	One can check that maps in the above collection glue together into a simplicial map which is an inverse of $f$.

\end{proof}

 \bibliographystyle{amsalpha}
\bibliography{GrConsQCat}

\providecommand{\bysame}{\leavevmode\hbox to3em{\hrulefill}\thinspace}
\providecommand{\MR}{\relax\ifhmode\unskip\space\fi MR }
\providecommand{\MRhref}[2]{%
  \href{http://www.ams.org/mathscinet-getitem?mr=#1}{#2}
}
\providecommand{\href}[2]{#2}
\begin{thebibliography}{DHMS04}

\bibitem[DHMS04]{DHKS}
William~G. Dwyer, P.~S. Hirschhorn, Kan~D. M, and J.~H. Smith, \emph{Homotopy
  {L}imit {F}unctors on {M}odel {C}ategories and {H}omotopical {C}ategories},
  2004.

\bibitem[Dwy95]{DS95}
J.~Dwyer, W. G.and~Spalinski, \emph{Homotopy theories and model categories}, in
  Handbook of Algebraic Topology, 1995.

\bibitem[GJ99]{GJ}
P.~G. Goerss and J.~F. Jardine, \emph{Simplicial {H}omotopy {T}heory},
  Birkhauser Verlag, 1999.

\bibitem[GZ67]{GZ}
P.~Gabriel and M.~Zisman, \emph{Calculus of {F}ractions and {H}omotopy
  {T}heory}, Ergebnisse der Mathematik imd ihrer Grenzgebiete, vol.~35,
  Springer-Verlag, 1967.

\bibitem[Hir02]{Hirchhorn}
Phillip~S. Hirchhorn, \emph{Model {C}ategories and their {L}ocalizations},
  Mathematical Surveys and Monographs, vol.~99, Amer. Math. Soc., Providence,
  RI, 2002.

\bibitem[HM15]{HM}
Gijs Heuts and Ieke Moerdijk, \emph{Left fibrations and homotopy colimits},
  Mathematische Zeitschrift \textbf{279} (2015), no.~3, 723--744.

\bibitem[Joh]{PJ}
P.~Johnston, \emph{Sketches of an elephant - {A} topos theory compendium}, vol.
  1 {\&} 2, Oxford University Press.

\bibitem[Joy08a]{AJ2}
A.~Joyal, \emph{Notes on quasi-categories},
  \url{http://www.math.uchicago.edu/~may/IMA/Joyal.pdf}, 2008.

\bibitem[Joy08b]{AJ1}
\bysame, \emph{Theory of quasi-categories and applications},
  \url{http://mat.uab.cat/~kock/crm/hocat/advanced-course/Quadern45-2.pdf},
  2008.

\bibitem[Lur09]{JL}
Jacob Lurie, \emph{Higher {T}opos {T}heory}, Annals of Mathematics Studies,
  vol. 170, Princeton University Press, Princeton, NJ, 2009.

\bibitem[Mey84]{Meyer}
J.~P. Meyer, \emph{Bar and cobar constructions {I}}, J. Pure and Applied
  Algebra \textbf{33} (1984), 163--207.

\bibitem[SGA61]{SGA}
\emph{Cat{\'e}gories fibr{\'e}es et descente}, Seminaire de g{\'e}ometrie
  alg{\'e}brique de l{'}Institut des Hautes {\'E}tudes Scientifiques (SGA 1),
  Paris, 1961.

\bibitem[Sha20]{sharma}
A.~Sharma, \emph{Symmetric monoidal categories and {$\Gamma$}-categories}, Th.
  and Appl. of categories \textbf{35} (2020), no.~14, 417--512.

\bibitem[Shu06]{Shul}
M.~Shulman, \emph{Homotopy limits and colimits and enriched homotopy theory},
  arXiv:math/0610194v3, 2006.

\bibitem[Tho79]{T79}
R.~W. Thomason, \emph{Homotopy colimits in the category of small categories},
  Mathematical Proceedings of the Cambridge Philosophical Society \textbf{85}
  (1979), no.~1, 91?109.

\end{thebibliography}

\end{document}